\documentclass{article}

\usepackage[french,british]{babel}

\usepackage[T1]{fontenc}

\usepackage{amsmath}

\usepackage{amsfonts}
\usepackage{amsbsy}
\usepackage{amssymb}

\usepackage{amsthm}

\theoremstyle{plain}

\usepackage{latexsym}
\usepackage{amsmath}
\usepackage{amsfonts}
\usepackage{amssymb}
\usepackage{psfrag}
\usepackage{graphicx}

\usepackage{amsmath,amssymb}
\usepackage{graphicx}

\usepackage{amsmath,amssymb}
\usepackage{graphicx}

\newtheorem{ozn}{Definition}[section]
\newtheorem{thm}{Theorem}[section]

\newtheorem{nas}{Corollary}[section]
\newtheorem{zau}{Remark}[section]
\newtheorem{lema}{Lemma}[section]

\newtheorem{remark}{Remark}[section]

\newcommand{\me}{\mathbf}
\newcommand{\mr}{\mathbb}
\newcommand{\mt}{\mathsf}
\newcommand{\md}{\mathcal}
\newcommand{\ld}{\left}
\newcommand{\rd}{\right}

\newcommand{\ip}{\int_{-\pi}^{\pi}}

\newcommand{\be}{\begin{equation}}
\newcommand{\ee}{\end{equation}}

\newcommand{\bem}{\begin{multline}}
\newcommand{\eem}{\end{multline}}

\newcommand{\bml}{\begin{multline*}}
\newcommand{\eml}{\end{multline*}}

\newcommand{\beg}{\begin{gather}}
\newcommand{\eeg}{\end{gather}}

\begin{document}

\title{Prediction problem for continuous time stochastic processes with periodically correlated increments observed with noise}

\author{
Maksym Luz\thanks {BNP Paribas Cardif in Ukraine, Kyiv, Ukraine, maksym.luz@gmail.com},
Mikhail Moklyachuk\thanks
{Department of Probability Theory, Statistics and Actuarial
Mathematics, Taras Shevchenko National University of Kyiv, Kyiv 01601, Ukraine, moklyachuk@gmail.com}
}

\date{\today}

\maketitle

\renewcommand{\abstractname}{Abstract}
\begin{abstract}
We propose solution of the problem of the mean square optimal estimation of linear functionals which depend on the unobserved values of
a continuous time  stochastic process with periodically correlated increments based on  observations of this process with periodically stationary noise.
To solve the problem, we transform the processes  to the sequences of stochastic functions  which form  an infinite dimensional vector stationary sequences.
In the case of known spectral densities of these sequences, we obtain formulas for calculating values of the mean square errors and the spectral characteristics of the optimal estimates of the functionals.
Formulas determining  the least favorable spectral densities and the minimax (robust) spectral
characteristics of the optimal linear estimates of functionals
are derived in the case where the  sets of admissible spectral densities are given.
\end{abstract}

\maketitle


\maketitle

\vspace{2ex}

\section*{Introduction}

In this article, we study the prediction problem for a continuous time stochastic process $\xi(t)$, $t\in\mr R$, with periodically
correlated increments $\xi^{(d)}(t,\tau T)=\Delta_{T\tau}^d \xi(t)$ of order $d$ and period $T$, where $\Delta_{s} \xi(t)=\xi(t)-\xi(t-s)$, based on observations of the process $\xi(t)$ with
a periodically correlated  noise stochastic process $\eta(t)$, $t\in\mr R$.
The resent studies, for example, by Basawa et al. \cite{Basawa}, Dudek et al. \cite{Dudek-Hurd}, Reisen et al. \cite{Reisen2018}, Solci et al. \cite{Solci}, show a  constant interest to the non-stationary models and robust methods of estimation.

Kolmogorov \cite{Kolmogorov},  Wiener \cite{bookWien} and  Yaglom \cite{Yaglom} developed effective methods of solution  of interpolation, extrapolation (prediction) and filtering problems for stationary stochastic sequences and processes. For a particular problem, they developed methods of finding an estimate $\widetilde{x}(t)$ constructed from available observations that minimizes the mean square error $\Delta(\widetilde{x}(t),f)=\mt E|x(t)-\widetilde{x}(t)|^2$ in the case where the spectral density $f(\lambda)$ of the stationary process or sequence $x(t)$ is exactly known and fixed.
Such estimates  are called optimal linear estimates  within this article.

The developed  classical estimation methods are not directly applicable in practice since the exact  spectral structure of the processes is not usually available. In this case the estimated spectral densities can be considered as the true ones. However, Vastola and Poor \cite{VastPoor1983} showed with the help of the concrete examples, that such substitution  can result in a significant increase of  the   estimate error.  Therefore it is reasonable to consider the estimates, called minimax-robust, which minimize the maximum of the mean-square errors for all spectral densities from a given set of admissible spectral densities simultaneously. The minimax-robust method of extrapolation was proposed by Grenander \cite{Grenander} who considered the estimation of the functional $Ax = \int_{0}^1 a (t) x (t) dt$ as a  game between two players, one of which minimizes the mean square error $\Delta(f,\widetilde{A}x)$ by $\widetilde{A}\zeta$ and another one maximizes the error by $f$. He showed that the game has a saddle point solution under proper conditions.
  For more details see the further study by Franke and Poor \cite{Franke1984}, Hosoya \cite{Hosoya} and  the survey paper by Kassam and Poor \cite{Kassam_Poor1985}.
  A wide range of results has been obtained by  Moklyachuk \cite{Moklyachuk:1981,Moklyachuk1981tv,Moklyachuk1982,Moklyachuk:2008,Moklyachuk2015}.
  These results have been  extended on the vector-valued stationary processes and sequences by Moklyachuk, Sidei and Masyutka \cite{Mokl_Mas_book,Sidei_book}.

The concept of stationarity   admits some generalizations, a combination of two of which --   stationary $d$th increments and periodical  correlation -- is in scope of this article.
Random processes with stationary $d$th increments $x(t)$ were introduced   by Yaglom and Pinsker \cite{Pinsker}. The increment sequence $x^{(d)}(t,\tau )=\Delta_{\tau}^d x(t)$ generated by such process is stationary by the variable $t$, namely, the mathematical expectations $\mt E\xi^{(d)}(t,\tau)$ and $\mt E\xi^{(d)}(t+s,\tau_1)\xi^{(d)}(t,\tau_2)$ do not depend on $t$.
Yaglom and Pinsker \cite{Pinsker} described the spectral representation of such process and the spectral density canonical factorization, and they also solved the extrapolation problem for these processes.
 The minimax-robust extrapolation, interpolation and filtering  problems for stochastic  processes with stationary increments were investigated by Luz and Moklyachuk \cite{Luz_Mokl_book}.

Dubovetska and Moklyachuk \cite{Dubovetska6}  derived the  classical and minimax-robust estimates for another generalization of stationary processes --  periodically  correlated (cyclostationary) processes, introduced by Gladyshev \cite{Glad1963}. The correlation function  $K(t,s)={\textsf{E}}{x(t)\overline{x(s)}}$ of such processes is a $T$-periodic function: $K(t,s)=K(t+T,s+T)$, which implies a time-dependent spectrum. Periodically  correlated  processes are widely used in signal processing and communications (see Napolitano \cite{Napolitano} for a review of recent works on cyclostationarity and its applications).

In this article, we deal with the problem of the mean-square optimal estimation of the linear functionals $A\xi=\int_{0}^{\infty}a(t)\xi(t)dt$ and $A_{NT}{\xi}=\int_{0}^{(N+1)T}a(t)\xi(t)dt$
which depend on the unobserved values of a continuous time stochastic process $\xi(t)$ with periodically stationary $d$-th
increments based on  observations of this process with periodically stationary noise at points $t<0$.

Similar problems for discrete time processes have been studied by Kozak, Luz and Moklyachuk \cite{Kozak_Mokl}, Luz and Moklyachuk \cite{Luz_Mokl_extra,Luz_Mokl_extra_GMI,Luz_Mokl_extra_noise_PCI}.
The problem of estimation of continuous time stochastic process $\xi(t)$ with periodically stationary $d$-th increments based on  observations of the process without noise at points $t<0$
was studied by Luz and Moklyachuk \cite{Luz_Mokl_extra_cont_PCI}.

The article is organized as follows.
   In section \ref{section_PC_process}, we describe a presentation of a continuous time periodically stationary process as a  stationary H-valued sequence.
   This approach is extended on the periodically stationary increments in the section \ref{section_PCI_c}.
  The traditional Hilbert space projection method of prediction is developed in section \ref{classical_extrapolation_e_n_c}.
    Particularly, formulas for calculating  the mean-square errors and the spectral characteristics of the optimal linear
estimates of the functionals $A\xi$ and $A_{NT}{\xi}$ are derived under some conditions on spectral densities.
An approach to solution of the prediction problem which is based on factorizations of spectral densities is developed
in section \ref{Prediction-factorizations}.
In section \ref{minimax_extrapolation} we present our results on minimax-robust prediction for the studied processes: relations that determine the least favourable spectral densities and the minimax spectral characteristics are derived for some classes of spectral densities.

\section{Continuous time periodically correlated processes and generated vector stationary sequences}
\label{section_PC_process}

In this section, we present a brief review of properties of periodically correlated processes and describe an approach to presenting it as stationary $H$-valued   sequences.
In the next section, this approach is applied  to develop the spectral theory for periodically correlated increment processes.

\begin{ozn}[Gladyshev \cite{Glad1963}] A mean-square continuous stochastic process
$\eta:\mathbb R\to H=L_2(\Omega,\mathcal F,\mathbb P)$, with $\textsf {E} \eta(t)=0,$ is called periodically correlated (PC) with period  $T$, if its correlation function  $K(t,s)={\textsf{E}}{\eta(t)\overline{\eta(s)}}$  for all  $t,s\in\mathbb R$ and some fixed $T>0$ is such that
\[
K(t,s)={\textsf{E}}{\eta(t)\overline{\eta(s)}}={\textsf{E}}{\eta(t+T)\overline{\eta(s+T)}}=K(t+T,s+T).
\]
\end{ozn}

For a periodically correlated  stochastic process $\eta(t)$, one can construct
the following sequence of stochastic functions \cite{DubovetskaMoklyachuk2013}, \cite{MoklyachukGolichenko2016}
\begin{equation} \label{zj}
\{\eta_j(u)=\eta(u+jT),u\in [0,T), j\in\mathbb Z\}.
\end{equation}
Sequence (\ref{zj}) forms a $L_2([0,T);H)$-valued
stationary sequence $\{\eta_j,j\in\mathbb Z\}$ with the correlation function
\begin{multline*}
B_{\eta}(l,j) = \langle\eta_l,\eta_j\rangle_H =\int_0^T
\textsf{E}[\eta(u+lT)\overline{\eta(u+jT)}]du
\\
=\int_0^TK_{\eta}(u+(l-j)T,u)du =
 B_{\eta}(l-j),
\end{multline*}
where
$K_{\eta}(t,s)=\textsf{E}{\eta}(t)\overline{{\eta}(s)}$
is the correlation function of the PC process $\eta(t)$.
Chose the following orthonormal basis in the space $L_2([0,T);\mathbb{R})$
\be\label{orthonormal_basis}
\{\widetilde{e}_k=\frac{1}{\sqrt{T}}e^{2\pi
i\{(-1)^k\left[\frac{k}{2}\right]\}u/T}, k=1,2,3,\dots\}, \quad
\langle \widetilde{e}_j,\widetilde{e}_k\rangle=\delta_{kj}.
\ee
Making use of this basis the stationary sequence $\{\eta_j,j\in\mathbb Z\}$  can be represented in the form
\be\label{eta_cont}
\eta_j= \sum_{k=1}^\infty \eta_{kj}\widetilde{e}_k,
\ee
where
\[\eta_{kj}=\langle\eta_j,\widetilde{e}_k\rangle =
\frac{1}{\sqrt{T}} \int_0^T \eta_j(v)e^{-2\pi
i\{(-1)^k\left[\frac{k}{2}\right]\}v/T}dv.\]
The sequence
$\{\eta_j,j\in\mathbb Z\},$
   or the corresponding to it  vector sequence
  \[\{\vec \eta_j=(\eta_{kj}, k=1,2,\dots)^{\top},
j\in\mathbb Z\},\]
is called a generated by the process $\{\eta(t),t\in\mathbb R\}$ vector stationary sequence.
The components $\{\eta_{kj}\}: k=1,2,\dots;j\in\mathbb Z$  of the generated stationary sequence
$\{\eta_j,j\in\mathbb Z\}$ satisfy the  relations \cite{Kallianpur}, \cite{Moklyachuk:1981}
\[
\textsf{E}{\eta_{kj}}=0, \quad \|{\eta}_j\|^2_H=\sum_{k=1}^\infty
\textsf{E}|\eta_{kj}|^2\leq P_\eta=B_\eta(0), \quad
\textsf{E}\eta_{kl}\overline{\eta}_{nj}=\langle
R_{\eta}(l-j)\widetilde{e}_k,\widetilde{e}_n\rangle.
\]
The correlation function $R_{\eta}(j)$  of the generated stationary
sequence $\{\eta_j,j\in\mathbb Z\}$
 is a correlation operator function.
 The correlation operator $R_{\eta}(0)=R_{\eta}$ is a
kernel operator and its kernel norm satisfies the following properties:
\[
\|{\eta}_j\|^2_H=\sum_{k=1}^\infty \langle R_{\eta}
\widetilde{e}_k,\widetilde{e}_k\rangle\leq P_\eta. \]
The generated stationary sequence $\{\eta_j,j\in\mathbb Z\}$ has the spectral density function
$g(\lambda)=\{g_{kn}(\lambda)\}_{k,n=1}^\infty$, that is positive valued operator  function of variable
 $\lambda\in [-\pi,\pi)$, if its correlation function $R_{\eta}(j)$ can be represented in the form
\[
\langle R_{\eta}(j)\widetilde{e}_k,\widetilde{e}_n\rangle=\frac{1}{2\pi} \int _{-\pi}^{\pi}
e^{ij\lambda}\langle g(\lambda) \widetilde{e}_k,\widetilde{e}_n\rangle d\lambda.
\]
We finish our review by the statement, that
for almost all   $\lambda\in [-\pi,\pi)$ the spectral density $f(\lambda)$ is a kernel operator with an integrable kernel norm
\[
\sum_{k=1}^\infty \frac{1}{2\pi} \int _{-\pi}^\pi \langle g(\lambda)
\widetilde{e}_k,\widetilde{e}_k\rangle d\lambda=\sum_{k=1}^\infty\langle R_{\eta}
\widetilde{e}_k,\widetilde{e}_k\rangle=\|{\eta}_j\|^2_H\leq P_\eta.
\]

\section{Stochastic processes with periodically correlated $d$th increments}
\label{section_PCI_c}

For a given stochastic process $\{\xi(t),t\in\mathbb R\}$, consider the stochastic $d$th increment process
\begin{equation}
\label{oznachPryrostu_cont}
\xi^{(d)}(t,\tau)=(1-B_{\tau})^d\xi(t)=\sum_{l=0}^d(-1)^l{d \choose l}\xi(t-l\tau),
\end{equation}
with the step $\tau\in\mr R$, generated by the stochastic process $\xi(t)$. Here $B_{\tau}$ is the backward shift operator: $B_{\tau}\xi(t)=\xi(t-\tau)$, $\tau\in \mr R$.

We prefer to use the notation $\xi^{(d)}(t,\tau)$ instead of widely used $\Delta_{\tau}^{d}\xi(t)$ to avoid a duplicate with the mean square error notation.

\begin{ozn}
\label{OznPeriodProc2_cont}
A stochastic process $\{\xi(t),t\in\mathbb R\}$ is called a  stochastic
process with periodically stationary (periodically correlated) increments  with the step $\tau\in\mr Z$ and the period $T>0$ if the mathematical expectations
\begin{eqnarray*}
\mt E\xi^{(d)}(t+T,\tau T) & = & \mt E\xi^{(d)}(t,\tau T)=c^{(d)}(t,\tau T),
\\
\mt E\xi^{(d)}(t+T,\tau_1 T)\xi^{(d)}(s+T,\tau_2 T)
& = & D^{(d)}(t+T,s+T;\tau_1T,v_2T)
\\
& = & D^{(d)}(t,s;\tau_1T,\tau_2T)
\end{eqnarray*}
exist for every  $t,s\in \mr R$, $\tau_1,\tau_2 \in \mr Z$ and for some fixed $T>0$.
\end{ozn}

The functions $c^{(d)}(t,\tau T)$ and  $D^{(d)}(t,s;\tau_1T,\tau_2 T)$ from the Definition \ref{OznPeriodProc2_cont} are called the \emph{mean value} and  the \emph{structural function} of the stochastic
process $\xi(t)$ with periodically stationary (periodically correlated) increments.

\begin{zau}
For spectral properties of one-pattern increment sequence $\chi_{\mu,1}^{(n)}(\xi(m)):=\xi^{(n)}(m,\mu)=(1-B_{\mu})^n\xi(m)$
see, e.g., \cite{Luz_Mokl_book}, pp. 1-8; \cite{Gihman_Skorohod}, pp. 48--60, 261--268; \cite{Yaglom}, pp. 390--430.
The corresponding results for continuous time increment process $\xi^{(n)}(t,\tau)=(1-B_{\tau})^n\xi(t)$ are described in \cite{Yaglom:1955}, \cite{Yaglom}.

For a review of the properties of periodically correlated processes, we refer to \cite{DubovetskaMoklyachuk2013}, \cite{MoklyachukGolichenko2016}. Here we present a generalization of these properties on the processes with periodically correlated increments from Definition \ref{OznPeriodProc2_cont}.
\end{zau}

For the stochastic process $\{\xi(t), t\in \mathbb R\}$ with periodically correlated  increments $\xi^{(d)}(t,\tau T)$ and the integer step $\tau$, we follow the procedure described in the Section \ref{section_PC_process} and construct
a sequence of stochastic functions
\begin{equation} \label{xj}
\{\xi^{(d)}_j(u):=\xi^{(d)}_{j,\tau}(u)=\xi^{(d)}_j(u+jT,\tau T),\,\, u\in [0,T), j\in\mathbb Z\}.
\end{equation}

Sequence (\ref{xj}) forms a $L_2([0,T);H)$-valued
stationary increment sequence $\{\xi^{(d)}_j,j\in\mathbb Z\}$ with the structural function
\begin{eqnarray*}
B_{\xi^{(d)}}(l,j)&= &\langle\xi^{(d)}_l,\xi^{(d)}_j\rangle_H =\int_0^T
\textsf{E}[\xi^{(d)}_j(u+lT,\tau_1 T)\overline{\xi^{(d)}_j(u+jT,\tau_2 T)}]du
\\&=&\int_0^T D^{(d)}(u+(l-j)T,u;\tau_1T,\tau_2T) du =
 B_{\xi^{(d)}}(l-j).
 \end{eqnarray*}
Making use of the orthonormal basis \eqref{orthonormal_basis} the stationary increment sequence  $\{\xi^{(d)}_j,j\in\mathbb Z\}$ can be represented in the form
\begin{equation} \label{zeta}
\xi^{(d)}_j= \sum_{k=1}^\infty \xi^{(d)}_{kj}\widetilde{e}_k,\end{equation}
where
\[\xi^{(d)}_{kj}=\langle\xi^{(d)}_j,\widetilde{e}_k\rangle =
\frac{1}{\sqrt{T}} \int_0^T \xi^{(d)}_j(v)e^{-2\pi
i\{(-1)^k\left[\frac{k}{2}\right]\}v/T}dv.\]

We call this sequence
$\{\xi^{(d)}_j,j\in\mathbb Z\},$
   or the corresponding to it vector sequence
\be \label{generted_incr_seq_i_n_c}
  \{\vec\xi^{(d)}(j,\tau)=\vec\xi^{(d)}_j=(\xi^{(d)}_{kj}, k=1,2,\dots)^{\top}=(\xi^{(d)}_{k}(j,\tau), k=1,2,\dots)^{\top},
j\in\mathbb Z\},\ee
\emph{an infinite dimension vector stationary increment sequence} generated by the increment process $\{\xi^{(d)}(t,\tau T),t\in\mathbb R\}$.
Further, we will omit the word vector in the notion generated vector stationary increment sequence.

Components $\{\xi^{(d)}_{kj}\}: k=1,2,\dots;j\in\mathbb Z$  of the generated stationary increment sequence
$\{\xi^{(d)}_j,j\in\mathbb Z\}$ are such that, \cite{Kallianpur}, \cite{Moklyachuk:1981}
\[
\textsf{E}{\xi^{(d)}_{kj}}=0, \quad \|\xi^{(d)}_j\|^2_H=\sum_{k=1}^\infty
\textsf{E}|\xi^{(d)}_{kj}|^2\leq P_{\xi^{(d)}}=B_{\xi^{(d)}}(0),
\]
and
\[\textsf{E}\xi^{(d)}_{kl}\overline{\xi^{(d)}_{nj}}=\langle
R_{\xi^{(d)}}(l-j;\tau_1,\tau_2)\widetilde{e}_k,\widetilde{e}_n\rangle.
\]
The \emph{structural function}  $R_{\xi^{(d)}}(j):=R_{\xi^{(d)}}(j;\tau_1,\tau_2)$  of the generated stationary increment
sequence $\{\xi^{(d)}_j,j\in\mathbb Z\}$
 is a correlation operator function.
 The correlation operator $R_{\xi^{(d)}}(0)=R_{\xi^{(d)}}$ is a
kernel operator and its kernel norm satisfies the following limitations:
\[
\|\xi^{(d)}_j\|^2_H=\sum_{k=1}^\infty \langle R_{\xi^{(d)}}
\widetilde{e}_k,\widetilde{e}_k\rangle\leq P_{\xi^{(d)}}. \]

 Suppose that  the structural function $R_{\xi^{(d)}}(j)$ admits a representation
\[
\langle R_{\xi^{(d)}}(j;\tau_1, \tau_2)\widetilde{e}_k,\widetilde{e}_n\rangle=\frac{1}{2\pi} \int _{-\pi}^{\pi}
e^{ij\lambda}(1-e^{-i\tau_1\lambda})^d(1-e^{i\tau_2\lambda})^d\frac{1}
{\lambda^{2d}}\langle f(\lambda) \widetilde{e}_k,\widetilde{e}_n\rangle d\lambda.
\]
Then
$f(\lambda)=\{f_{kn}(\lambda)\}_{k,n=1}^\infty$ is a \emph{spectral density function} of the generated stationary increment sequence $\{\xi^{(d)}_j,j\in\mathbb Z\}$. It  is a positive valued operator  functions of variable
 $\lambda\in [-\pi,\pi)$, and for almost all   $\lambda\in [-\pi,\pi)$ it is a kernel operator with an integrable kernel norm
\begin{multline} \label{P1}
\sum_{k=1}^\infty \frac{1}{2\pi} \int _{-\pi}^\pi (1-e^{-i\tau_1\lambda})^d(1-e^{i\tau_2\lambda})^d\frac{1}
{\lambda^{2d}} \langle f(\lambda)
\widetilde{e}_k,\widetilde{e}_k\rangle d\lambda
\\
=
\sum_{k=1}^\infty\langle R_{\xi^{(d)}}
\widetilde{e}_k,\widetilde{e}_k\rangle=\|{\zeta}_j\|^2_H\leq P_{\xi^{(d)}}.
\end{multline}

The stationary $d$th increment sequence $\vec{\xi}^{(d)}_j$ admits the spectral representation \cite{Karhunen}
\[
\vec{\xi}^{(d)}_j=\int_{-\pi}^{\pi}e^{i \lambda j}(1-e^{-i\tau\lambda})^{d}\frac{1}{(i\lambda)^{d}}d\vec{Z}_{\xi^{(d)}}(\lambda),
\]
where $\vec{Z}_{\xi^{(d)}}(\lambda)=\{Z_{k}(\lambda)\}_{k=1}^{\infty}$ is a vector-valued random process with uncorrelated increments on $[-\pi,\pi)$.

Consider the generated stationary stochastic sequence $\vec\eta_j$ defined in Section \ref{section_PC_process_i_n_c}, which is uncorrelated with the increment sequence $\vec \xi^{(d)}_j$. It admits the spectral representation
\begin{equation}
\label{SpectrPred_incr_eta_c}
\vec\eta_j=\int_{-\pi}^{\pi}e^{i\lambda j}d\vec{Z}_{\eta}(\lambda),
\end{equation}
where $\vec{Z}_{\eta}(\lambda)$   is a vector-valued random process with uncorrelated increments on $[-\pi,\pi)$.
The spectral representation  of the
sequence $\vec \zeta^{(d)}_j$, generated by the process $\zeta(t)=\xi(t)+\eta(t)$, is determined by the spectral densities $f(\lambda)$
and $g(\lambda)$ by the relation
\begin{equation}
\label{SpectrPred_incr_zeta_c}
\vec{\zeta}^{(d)}_j=\int_{-\pi}^{\pi}e^{i \lambda j}(1-e^{-i\tau\lambda})^{d}\frac{1}{(i\lambda)^{d}}d\vec{Z}_{\xi^{(d)}+\eta^{(d)}}(\lambda).
\end{equation}
The random processes $\vec{Z}_{\eta}(\lambda)$ and $\vec{Z}_{\eta^{(d)}}(\lambda)$ are connected by the relation $d\vec{Z}_{\eta^{(d)}}(\lambda)=(i\lambda)^d d\vec{Z}_{\eta }(\lambda)$,
$\lambda\in[-\pi,\pi)$, see \cite{Luz_Mokl_book}. The spectral density $p(\lambda)=\{p_{kn}(\lambda)\}_{k,n=1}^\infty$ of the
sequence $\vec \zeta^{(d)}_j$  is determined by the spectral densities $f(\lambda)$
and $g(\lambda)$ by the relation
\[p(\lambda)=f(\lambda)+\lambda^{2d}g(\lambda).\]

In the space $ H=L_2(\Omega, \cal F, \mt P)$, consider    a closed linear subspace
\[
H(\vec{\xi}^{\,(d)})=\overline{span}\{\xi^{(d)}_{kj}: k=1,2,\dots;\, j\in \mathbb Z \}\]
 generated by the components
of the generated stationary increment sequence
$\vec{\xi}^{\,(d)}=\{\xi^{(d)}_{kj}=\xi^{(d)}_{k}(j,\tau),\,\tau>0\}$.
For $q\in \mathbb Z$, consider also a closed linear subspace
\[H^{q}(\vec{\xi}^{\,(d)})=\overline{span}\{\xi^{(d)}_{kj}: k=1,2,\dots;\, j\leq q \}.\]
Define a subspace
\[ S(\vec{\xi}^{\,(d)})=\bigcap_{q\in \mathbb{Z}} H^{q}(\vec{\xi}^{\,(d)})
\]
  of the Hilbert space $H(\vec{\xi}^{\,(d)})$. The  space $H(\vec{\xi}^{\,(d)})$ admits a decomposition
$ H(\vec{\xi}^{\,(d)})=S(\vec{\xi}^{\,(d)})\oplus R(\vec{\xi}^{\,(d)}) $
where $R(\vec{\xi}^{\,(d)})$ is the orthogonal complement of the subspace $S(\vec{\xi}^{\,(d)})$ in the space $H(\vec{\xi}^{\,(d)})$.
\begin{ozn}
A stationary (wide sense)    increment sequence $\vec\xi^{(d)}_{j}=\{\xi^{(d)}_{kj}\}_{k=1}^{\infty}$ is called regular if $H(\vec{\xi}^{\,(d)})=R(\vec{\xi}^{\,(d)})$,
and it is called singular if
$H(\vec{\xi}^{\,(d)})=S(\vec{\xi}^{\,(d)})$.
\end{ozn}

\begin{thm}
A stationary    increment sequence $\xi^{(d)}_{j}$ is uniquely represented in the form
\begin{equation} \label{rozklad_cont}
\xi^{(d)}_{kj}
=\xi^{(d)}_{S,kj}+\xi^{(d)}_{R,kj}
\end{equation}
where  $\xi^{(d)}_{R,kj}, k=1,\ldots,\infty$, is a regular stationary     increment sequence and
$\xi^{(d)}_{S,kj}, k=1,\dots,\infty$, is a singular stationary   increment sequence.
The   increment sequences
$\xi^{(d)}_{R,kj}$  and
$\xi^{(d)}_{S,kj}$
are
orthogonal for all $ j\in\mathbb{Z} $. They are defined by the formulas
\begin{align*} \xi^{(d)}_{S,kj}&=\mt E[\xi^{(d)}_{kj}|S(\vec{\xi}^{\,(d)})],
\\
\xi^{(d)}_{R,kj}&=\xi^{(d)}_{kj}
-\xi^{(d)}_{S,kj}.
\end{align*}
\end{thm}

Consider an innovation sequence  ${\vec\varepsilon(u)=\{\varepsilon_m(u)\}_{m=1}^M, u \in\mathbb Z}$ for a regular stationary
increment, namely, the sequence of uncorrelated random
variables such that $\mathsf{E} \varepsilon_m(u)\overline{\varepsilon}_j(v)=\delta_{mj}\delta_{uv}$,   $\mathsf{E} |\varepsilon_m(u)|^2=1, m,j=1,\dots,M; u \in\mathbb Z$, and $H^{r}(\vec\xi^{(d)} )=H^{r}(\vec\varepsilon)$ holds true for all $r \in \mathbb Z$, where  $H^r(\vec\varepsilon)$ is
the Hilbert space generated by elements $ \{ \varepsilon_m(u):m=1,\dots,M; u\leq r\}$,
 $\delta_{mj}$ and $\delta_{uv}$ are Kronecker symbols.

\begin{thm}\label{thm 4_cont}
A   stationary   increment sequence
$\vec\xi^{(d)}_{j}$ is regular if and only if there exists an
innovation sequence ${\vec\varepsilon(u)=\{\varepsilon_m(u)\}_{m=1}^M, u \in\mathbb Z}$
and a sequence of matrix-valued
functions $\varphi^{(d)}(l,\tau) =\{\varphi^{(d)}_{km}(l,\tau) \}_{k=\overline{1,\infty}}^{m=\overline{1,M}}$, $l\geq0$, such that
\begin{equation}\label{odnostRuhSer_cont}
\sum_{l=0}^{\infty}
\sum_{k=1}^{\infty}
\sum_{m=1}^{M}
|\varphi^{(d)}_{km}(l,\tau)|^2
<\infty,\quad
\xi^{(d)}_{kj}=
\sum_{l=0}^{\infty}\sum_{m=0}^{M}\varphi^{(d)}_{km}(l,\tau)\vec\varepsilon_m(j-l).
\end{equation}
Representation (\ref{odnostRuhSer_cont}) is called the canonical
moving average representation of the generated stationary   increment
sequence $\vec\xi^{\,(d)}_{j}$.
\end{thm}

The spectral function $F(\lambda)$ of  a stationary    increment sequence $\vec\xi^{\,(d)}_{j}$
which admits   canonical representation
 $(\ref{odnostRuhSer_cont})$  has the spectral density  $f(\lambda)=\{f_{ij}(\lambda)\}_{i,j=1}^\infty$ admitting the canonical
factorization
\begin{equation}\label{SpectrRozclad_f_cont}
f(\lambda)=
\varphi(e^{-i\lambda})\varphi^*(e^{-i\lambda}),
\end{equation}
where the function
$\varphi(z)=\sum_{k=0}^{\infty}\varphi(k)z^k$ has
analytic in the unit circle $\{z:|z|\leq1\}$
components
$\varphi_{ij}(z)=\sum_{k=0}^{\infty}\varphi_{ij}(k)z^k; i=1,\dots,\infty; j=1,\dots,M$. Based on moving average  representation
$(\ref{odnostRuhSer_cont})$ define
\[\varphi_{\tau}(z)=\sum_{k=0}^{\infty} \varphi^{(d)}(k,\tau)z^k=\sum_{k=0}^{\infty}\varphi_{\tau}(k)z^k.\]
 Then the following factorization holds true:
\be
 \frac{|1-e^{i\lambda\tau}|^{2d}}{\lambda^{2d}}f(\lambda)= \varphi_{\tau}(e^{-i\lambda})
\varphi^*_{\tau}(e^{-i\lambda}),\quad \varphi_{\tau}(e^{-i\lambda})=\sum_{k=0}^{\infty}\varphi_{\tau}(k)e^{-i\lambda k}.
 \label{dd_cont}
\ee

\section{Hilbert space projection method of prediction}\label{classical_extrapolation_e_n_c}

Let a periodically correlated increment process $\xi^{(d)}(t,\tau T)$, $t\in\mathbb R$, generates by formula (\ref{zeta})
an infinite dimension vector stationary increment sequence $\{\vec{\xi}^{(d)}_j,j\in\mr Z\}$ which has
 the spectral density matrix $f(\lambda)=\{f_{ij}(\lambda)\}_{i,j=1}^{\infty}$.
As a noise process, consider a periodically stationary stochastic process $\eta(t)$, $t\in\mathbb R$, uncorrelated with
the process $\xi(t)$. Let the process $\eta(t)$  generates by formula (\ref{eta_cont}) an infinite dimension vector stationary
sequence $\{\vec{\eta}_j,j\in\mr Z\}$  with the spectral density matrix
$g(\lambda)=\{g_{ij}(\lambda)\}_{i,j=1}^{\infty}$.

By the classical \textbf{prediction problem} we understand the problem of the mean square optimal linear  estimation  of the functionals
\[A{\xi}=\int_{0}^{\infty}a(t)\xi(t)dt,\quad A_{NT}{\xi}=\int_{0}^{(N+1)T}a(t)\xi(t)dt\]
which depend on the unknown values of the stochastic process $\xi(t)$.
Estimates are based on observations of the process $\zeta(t)=\xi(t)+\eta(t)$  at points $t<0$.

Assumptions:
\begin{itemize}
  \item  the mean
values of the increment sequence $\vec{\xi}^{(d)}_j$ and stationary
sequence  $\vec{\eta}_j$ equal to 0; the increment step $\tau>0$;
  \item the spectral densities $f(\lambda)$ and $g(\lambda)$ satisfy the \emph{minimality condition}
\be
 \ip \text{Tr}\left[ \frac{\lambda^{2d}}{|1-e^{i\lambda\tau}|^{2d}}(f(\lambda)+\lambda^{2d}g(\lambda))^{-1}\right]
 d\lambda<\infty.
\label{umova11_e_n_c}
\ee
\end{itemize}

    The latter assumption is the necessary and sufficient condition under which the mean square errors of the optimal
    estimates of the functional $A\vec\xi$ to be defined below is not equal to $0$.

The  Hilbert space projection method of estimation may be applied under the condition that the  element, which we want to estimate,
  belongs to the Hilbert space   $ H=L_2(\Omega, \mathcal F, \mt P)$   of random variables  with a zero mean value and a finite variance.
  That is not the case for the functional $A{\xi}$.
  To overcome this difficulty we find a representation of
   the functional $A{\xi}$
    as a sum of a functional with finite second moment  from the space $H$ and a functional depended on the  observed values of the process $\zeta(t)=\xi(t)+\eta(t)$.
    This representation is described by the following two lammas.

\begin{lema}[\cite{Luz_Mokl_extra_cont_PCI}]
\label{predst A_cont}
The  linear functional
\[A\zeta=\int_{0}^{\infty}a(t)\zeta(t)dt\]
allows the representation
\[
    A\zeta=B\zeta-V\zeta,\]
where
\[
    B\zeta=\int_0^{\infty} b^{\tau}(t)\zeta^{(d)}(t,\tau T)dt,\quad
    V\zeta=\int_{-\tau T d}^{0}v^{\tau}(t)\zeta(t)dt,\]
and
\begin{eqnarray} \label{koefv_cont}
    v^{\tau}(t)&=&\sum_{l=\left\lceil-\frac{t}{\tau T}\right\rceil}^d(-1)^l{d \choose l}b^{\tau}(t+l\tau T),\quad t\in [-\tau T d;0),
\\
    \label{koef b_cont}b^{\tau}(t)&=&\sum_{k=0}^{\infty}a(t+\tau T k)d(k)=D^{\tau T}\me a(t) ,\,\,t\geq0,
\end{eqnarray}
Here  $\lceil x\rceil$ denotes the least integer greater than or equal to  $x$, $[x]$ denotes the integer part of $x$, coefficients
  $\{d(k):k\geq0\}$ are determined by the relation
\[\sum_{k=0}^{\infty}d (k)x^k=\left(\sum_{j=0}^{\infty}x^{
j}\right)^d,\]
$D^{\tau T}$ is the linear transformation acting on an arbitrary function $x(t)$, $t\geq0$, as follows:
\[
    D^{\tau T}\me x(t)=\sum_{k=0}^{\infty}x(t+\tau T k )d(k).\]
\end{lema}

From Lemma \ref{predst A_cont}, we obtain the following representation of the functional $A\xi$:
\[A\xi=A\zeta-A\eta=B\zeta-A\eta-V\zeta=H\xi-V\zeta,\]
where
\[
 H\xi:=B\zeta-A\eta,\]
and
\[
 A\zeta=\int_{0}^{\infty}a(t)\zeta(t)dt,\quad A\eta=\int_{0}^{\infty}a(t)\eta(t)dt,\]
\[
 B\zeta=\int_0^{\infty} b^{\tau}(t)\zeta^{(d)}(t,\tau T)dt,\quad
 V\zeta=\int_{-\tau T d}^{0}v^{\tau}(t)\zeta(t)dt, \]
the functions $b_{\tau}(t)$, $t\in[0;\infty)$, and $v_{\tau}(t)$, $t\in[-\tau
Td;0)$, are calculated by formulas (\ref{koef b_cont}) and
(\ref{koefv_cont}) respectively.

The functional $H\xi$ allows a representation  in terms of the sequences   $\vec\eta_j=(\eta_{kj},k=1,2,\dots)^{\top}$ and $\vec\zeta^{(d)}_j=\vec\xi^{(d)}_j+\vec\eta^{(d)}_j=(\zeta^{(d)}_{kj},k=1,2,\dots)^{\top}$, $j\in\mr Z$, which is described in the following lemma.

\begin{lema}\label{predst H_cont}
The functional $H\xi=B\zeta-A\eta$ can be represented in the form
\[
 H\xi= \sum_{j=0}^{\infty}
{(\vec{b}^{\tau}_j)}^{\top}\vec{\zeta}^{(d)}_j- \sum_{j=0}^{\infty}
{(\vec{a}_j)}^{\top}\vec{\eta}_j
=B\vec \zeta - V\vec \eta=:H\vec \xi,
\]
where the vector
\[
\vec{b}^{\tau}_j=(b^{\tau}_{kj},k=1,2,\dots)^{\top}=
(b^{\tau}_{1j},b^{\tau}_{3j},b^{\tau}_{2j},\dots,b^{\tau}_{2k+1,j},b^{\tau}_{2k,j},\dots)^{\top},
\]
with the entries
\[
 b^{\tau}_{kj}=\langle b^{\tau}_j,\widetilde{e}_k\rangle =
\frac{1}{\sqrt{T}} \int_0^T b^{\tau}_j(v)e^{-2\pi
i\{(-1)^k\left[\frac{k}{2}\right]\}v/T}dv,
\]
and the vector
\[
\vec{a}_j=(a_{kj},k=1,2,\dots)^{\top}=
(a_{1j},a_{3j},a_{2j},\dots,a_{2k+1,j},a_{2k,j},\dots)^{\top},
\]
with the entries
\begin{multline*}
a_{kj}=\langle a_j,\widetilde{e}_k\rangle =
\frac{1}{\sqrt{T}} \int_0^T a_j(v)e^{-2\pi
i\{(-1)^k\left[\frac{k}{2}\right]\}v/T}dv,
\\
 k=1,2,\dots,\,j=0,1,\ldots,\infty.
 \end{multline*}

The coefficients $\{\vec{a}_j, j=0,1,\dots,\infty\}$ and  $\{\vec{b}^{\tau}_j, j=0,1,\dots,\infty\}$ are related as
\be \label{a_b_relation_e_n_c}
\vec{b}_j^{\tau}=\sum_{m=j}^{\infty}\mt{diag}_{\infty}(d_{\tau}(m-j))\vec{a}_m=(D^{\tau}{\me a})_j,  \quad j=0,1,\dots,\infty.
\ee
where
$\me a=((\vec{a}_0)^{\top},(\vec{a}_1)^{\top}, \ldots)^{\top}$, the coefficients $\{d_{\tau}(k):k\geq0\}$ are determined by the relationship
\[
 \sum_{k=0}^{\infty}d_{\tau}(k)x^k=\left(\sum_{j=0}^{\infty}x^{\tau j}\right)^d,\]
$D^{\tau}$ is a linear transformation  determined by a matrix with the infinite dimension matrix entries
 $D^{\tau}(k,j), k,j=0,1,\dots$ such that $D^{\tau}(k,j)=\mt{diag}_{\infty}(d_{\tau}(j-k))$ if $0\leq k\leq j \leq \infty$ and $D^{\tau}(k,j)=\mt{diag}_{\infty}(0)$ for $0\leq j<k\leq \infty$; $\mt{diag}_{\infty}(x)$ denotes an infinite dimensional diagonal matrix with the entry $x$ on its diagonal.
\end{lema}
\begin{proof}
See Appendix.
\end{proof}

Assume, that coefficients $\{\vec{a}_j, j=0,1,\dots\}$ and   $\{\vec{b}^{\tau}_j, j=0,1,\dots\}$, that determine the functional $H\vec\xi$, satisfy the conditions
\begin{equation} \label{coeff_a_e_n_c}
 \sum_{j=0}^\infty\|
 \vec{a}_j\|<\infty, \quad
 \sum_{j=0}^\infty(j+1)\|
 \vec{a}_j\|<\infty, \quad
  \|\vec{a}_j\|^2=\sum_{k=1}^\infty |{a}_{kj}|^2,
\end{equation}
\begin{equation} \label{coeff_b_e_n_c}
 \sum_{j=0}^\infty\|
 \vec{b}^{\tau}_j\|<\infty, \quad
 \sum_{j=0}^\infty(j+1)\|
 \vec{b}^{\tau}_j\|<\infty, \quad
  \|\vec{b}^{\tau}_j\|^2=\sum_{k=1}^\infty |{b}^{\tau}_{kj}|^2.
\end{equation}

 Under conditions (\ref{coeff_a_e_n_c}) - (\ref{coeff_b_e_n_c}) the functional $H\vec{\xi}$ has finite second moment. Since the functional  $V\zeta$ depends on the  observations  $ \{\xi(t)+\eta(t):t\in[\tau T d;0)\}$, the estimates $\widehat{A}\xi$ and $\widehat{H}\vec\xi$ of the functionals $A\xi$ and $H\vec\xi$, as well as  the mean-square errors $\Delta(f,g;\widehat{A}\xi)=\mt E |A\xi-\widehat{A}\xi|^2$ and $\Delta(f,g;\widehat{H}\vec\xi)=\Delta(f,g;\widehat{H}\vec\xi)=\mt E
|H\vec\xi-\widehat{H}\vec\xi|^2$ of the estimates $\widehat{A}\xi$ and $\widehat{H}\vec\xi$ satisfy the  relations
\be\label{mainformula_e_n_c}
    \widehat{A}\xi=\widehat{H}\vec\xi-V\zeta,\ee
\begin{align*}
  \Delta(f,g;\widehat{A}\xi)
    &=\mt E |A\xi-\widehat{A}\xi|^2
    \\
    &=\mt E|H\vec\xi-V\zeta-\widehat{H}\vec\xi+V\zeta|^2
    \\
    &=\mt E|H\vec\xi-\widehat{H}\vec\xi|^2
    \\
    &=\Delta(f,g;\widehat{H}\vec\xi).
 \end{align*}
Thus, the  functional $H\vec\xi$ is a target element to be estimated. Let us describe its spectral representation. Making use of representations (\ref{SpectrPred_incr_zeta_c}) and (\ref{SpectrPred_incr_eta_c}), we  obtain
\[H\vec\xi=\int_{-\pi}^{\pi}(\vec{B}_{\tau}(e^{i\lambda}))^{\top}
\frac{(1-e^{-i\lambda \tau})^d}{(i\lambda)^d}d\vec{Z}_{\xi^{(d)}+\eta^{(d)}}(\lambda)
-\int_{-\pi}^{\pi}(\vec{A}(e^{i\lambda}))^{\top}d\vec{Z}_{\eta}(\lambda),\]
where
\[\vec{B}_{\tau}(e^{i\lambda})=\sum_{j=0}^{\infty}\vec{b}_j^{\tau}e^{i\lambda j}=\sum_{j=0}^{\infty}(D^{\tau}\me a)_je^{i\lambda j},\quad \vec{A}(e^{i\lambda})=\sum_{j=0}^{\infty}\vec{a}_je^{i\lambda j}.\]

The classical approach of estimation consists in finding a projection of the element $H\vec\xi$ on  the closed linear subspace of  $H=L_2(\Omega,\mathcal{F},\mt P)$ generated by the observations. Let us define this subspace as
\[H^{0-}(\xi^{(d)}_{\tau}+\eta^{(d)}_{\tau})
=\overline{span}\{\vec\xi^{(d)}_{kj}+\vec\eta^{(d)}_{kj}:k=1,\dots,\infty;j=-1,-2,-3,\dots\}.\]
Define also the    closed linear subspaces of the Hilbert space
$L_2(f(\lambda)+\lambda^{2d}g(\lambda))$  of vector-valued functions endowed by the inner product \[\langle g_1;g_2\rangle=\ip (g_1(\lambda))^{\top}(f(\lambda)+\lambda^{2d}g(\lambda))\overline{g_2(\lambda)}d\lambda\]
as
\begin{multline*}
L_2^{0-}(f(\lambda)+\lambda^{2d}g(\lambda))=\overline{span}\{
 e^{i\lambda j}(1-e^{-i\lambda \tau})^d\frac{1}{(i\lambda)^d}\vec{\delta}_{k},
 \\
  k=1,2,3,\dots;\,\, j=-1,-2,-3,\dots\},
 \end{multline*}
where $\vec{\delta}_k=\{\delta_{kl}\}_{l=1}^{\infty}$, $\delta_{kl}$ are Kronecker symbols.

\begin{remark}
  Representation
 (\ref{SpectrPred_incr_zeta_c})
yields a  map between the elements $e^{i\lambda j}(1-e^{-i\lambda\tau})^d(i\lambda)^{-d}\vec{\delta}_k$ of the space
$L_2^{0-}(f(\lambda)+\lambda^{2d}g(\lambda))$
and the elements
$\vec\xi^{(d)}_{kj}+\vec\eta^{(d)}_{kj}$
of the space
$H^{0-}(\xi^{(d)}_{\tau}+\eta^{(d)}_{\tau})$.
\end{remark}

 The mean square optimal estimate
$\widehat{H}\vec\xi$ is found  as a projection of the element $H\vec\xi$ on the
subspace $H^{0-}(\xi^{(d)}_{\tau}+\eta^{(d)}_{\tau})$:
\[
\widehat{H}\vec\xi=\mt{Proj}_{H^{0-}(\xi^{(d)}_{\tau}+\eta^{(d)}_{\tau})}H\vec\xi.
\]
Relation (\ref{mainformula_e_n_c}) let us write  the optimal estimate $\widehat{A}\xi$
in the form
\[
\widehat{A}\xi=\mt{Proj}_{H^{0-}(\xi^{(d)}_{\tau}+\eta^{(d)}_{\tau})}H\vec\xi-V\zeta
\]
or
in the form
\be \label{otsinka_A_e_n_c}
 \widehat{A}\xi=\ip
(\vec{h}_{\tau}(\lambda))^{\top}d\vec{Z}_{\xi^{(d)}+\eta^{(d)}}(\lambda)-\int_{-\tau T d}^{0}v^{\tau}(t)\zeta(t)dt,
\ee
 where
$\vec{h}_{\tau}(\lambda)=\{h_{k}^{\tau}(\lambda)\}_{k=1}^{\infty}$ is the spectral characteristic of the   estimate $\widehat{H}\vec\xi$.

Denote
\[
\vec{A}_{\tau}(e^{i\lambda})=(1-e^{i\lambda \tau})^d\vec{A}(e^{i\lambda})=\sum_{j=0}^{\infty}\vec{a}_j^{\tau}e^{i\lambda j},
\]
where
\be\label{coeff a_mu_e_n_c}
    \vec{a}_j^{\tau}
    =\sum_{l=0}^{\min\ld\{\ld[j / \tau\rd],d\rd\}}
    (-1)^l{d \choose l}\vec{a}(j-\tau l),\quad j\geq0.\ee
Define the vector
\[
     \me a^{\tau}=((\vec{a}^{\tau}_0)^{\top},(\vec{a}^{\tau}_1)^{\top},(\vec{a}^{\tau}_2)^{\top}, \ldots)^{\top}.\]
With the help of
  the Fourier coefficients
\begin{align*}
  T^{\tau}_{l,j} &  =\frac{1}{2\pi}\int_{-\pi}^{\pi}
e^{i\lambda(j-l)}
\dfrac{\lambda^{2d}}{|1-e^{i\lambda\tau}|^{2d}}\ld[g(\lambda)
(f(\lambda)+\lambda^{2d}g(\lambda))^{-1}\rd]^{\top}
d\lambda,\quad l,j\geq0,
\\
 P_{l,j}^{\tau} &=\frac{1}{2\pi}\int_{-\pi}^{\pi} e^{i\lambda (j-l)}
\dfrac{\lambda^{2d}}{|1-e^{i\lambda\tau}|^{2d}}\ld[
(f(\lambda)+\lambda^{2d}g(\lambda))^{-1}\rd]^{\top}
d\lambda,\quad l,j\geq0,
\\
   Q_{l,j}&=\frac{1}{2\pi}\int_{-\pi}^{\pi}
e^{i\lambda(j-l)}\ld[f(\lambda)(f(\lambda)+\lambda^{2d}g(\lambda))^{-1}g(\lambda)\rd]^{\top}d\lambda,\quad l,j\geq0,
\end{align*}
of the corresponding matrix functions, define
the linear operators $\me P_{\tau}$, $\me T_{\tau}$ and $\me Q$ in the space $\ell_2$  by matrices with the infinite dimensional matrix entries
 $(\me P_{\tau})_{l,j}=P_{l,j}^{\tau}$,   $(\me T_{\tau})_{l, j} =T^{\tau}_{l,j}$ and $(\me Q)_{l,j}=Q_{l,j}$, $l,j\geq0$.

Notation:     $\langle \vec x, \vec y\rangle=\sum_{j=0}^{\infty}(\vec x_j)^{\top}\overline{\vec y}_j$ for vectors $\vec x=((\vec x_0)^{\top},(\vec x_1)^{\top},(\vec x_2)^{\top},\ldots)^{\top}$, $\vec y=((\vec y_0)^{\top},(\vec y_1)^{\top},(\vec y_2)^{\top},\ldots)^{\top}$.

\begin{thm}\label{thm1_e_n_c}
Consider two uncorrelated
processes: a stochastic
process $\xi(t)$, $t\in \mr R$ with a periodically stationary  increments, which determines a
generated stationary  $d$th increment sequence
$\vec{\xi}^{(d)}_j$ with the spectral density matrix $f(\lambda)=\{f_{kn}(\lambda)\}_{k,n=1}^{\infty}$, and   a periodically stationary stochastic
process $\eta(t)$, $t\in \mr R$,  which determines
   a
generated stationary  sequence
$\vec{\eta}_j$ with the spectral density matrix $g(\lambda)=\{g_{kn}(\lambda)\}_{k,n=1}^{\infty}$.
Let the coefficients $\vec {a}_j$, $ \vec{b}^{\tau}_j$, $j=0,1,\ldots,$ generated by the function $a(t)$, $t\geq 0$, satisfy conditions  (\ref{coeff_a_e_n_c}) -- (\ref{coeff_b_e_n_c}).
 Let  minimality condition
(\ref{umova11_e_n_c}) be satisfied.
The optimal linear
estimate $\widehat{A}\xi$ of the functional $A\xi$ based on
observations of the process $\xi(t)+\eta(t)$ at points   $t<0$ is calculated by
formula (\ref{otsinka_A_e_n_c}). The spectral characteristic
   $\vec{h}_{\tau}(\lambda)$ is
calculated by formula
\begin{align}
\notag
(\vec{h}_{\tau}(\lambda))^{\top}&=(\vec{B}_{\tau}(e^{i\lambda}))^{\top}
\frac{(1-e^{-i\lambda \tau})^d}{(i\lambda)^d}
\\
&\quad-
\ld((\vec{A}_{\tau}(e^{i\lambda}))^{\top}g(\lambda)+(\vec{C}_{\tau}(e^{i\lambda}))^{\top}\rd)\frac{(-i\lambda)^{d}}{(1-e^{i\lambda \tau})^d}
(f(\lambda)+\lambda^{2d}g(\lambda))^{-1},
 \label{spectr_A_e_n_c}
 \end{align}
where
\[ \vec{C}_{\tau}(e^{i \lambda})=\sum_{j=0}^{\infty}
(\me P_{\tau}^{-1}D^{\tau}\me a-\me P_{\tau}^{-1}\me T_{\tau}\me a^{\tau})_j
e^{i\lambda j}.\]
The value of the
mean-square error is calculated by the
formula
\begin{align}
\notag \Delta(f,g;\widehat{A}\vec\xi)&=\frac{1}{2\pi}\int_{-\pi}^{\pi}
\frac{\lambda^{2d}}{|1-e^{i\lambda\tau}|^{2d}}
(\me C^{f}_{\tau}(e^{i\lambda}))^{\top}(f(\lambda)+{\lambda}^{2d}g(\lambda))^{-1}f(\lambda)
\\\notag
&\quad\quad\times
(f(\lambda)+{\lambda}^{2d}g(\lambda))^{-1}
\overline{\me C^{f}_{\tau}(e^{i\lambda})}
d\lambda
\\\notag
&\quad+\frac{1}{2\pi}\int_{-\pi}^{\pi}
\frac{\lambda^{4n}}{|1-e^{i\lambda\tau}|^{4n}}
(\me C^{g}_{\tau}(e^{i\lambda}))^{\top}(f(\lambda)+{\lambda}^{2d}g(\lambda))^{-1} g(\lambda)
\\\notag
&\quad\quad\times
(f(\lambda)+{\lambda}^{2d}g(\lambda))^{-1}
\overline{\me C^{g}_{\tau}(e^{i\lambda})}
d\lambda
\\
 &=
\langle D^{\tau}\me a- \me T_{\tau}\me
 a^{\tau},\me P_{\tau}^{-1}D^{\tau}\me a-\me P_{\tau}^{-1}\me T_{\tau}\me a^{\tau}\rangle+\langle\me Q\me a,\me
 a\rangle.\label{pohybkaA_e_n_c}
\end{align}
where
\begin{align*}
{\me C^{f}_{\tau}(e^{i\lambda})}
&:=
\overline{g(\lambda)}\vec{A}_{\tau}(e^{i\lambda}) +
\vec{C}_{\tau}(e^{i \lambda}),
\\
{\me C}^{g}_{\tau}(e^{i \lambda})
&:=
{|1-e^{i\lambda\tau}|^{2d}}\lambda^{-2d}\overline{f^0(\lambda)}\vec{A}(e^{i\lambda})
-{(1-e^{-i\lambda\tau})^d}
\vec{C}_{\tau}(e^{i \lambda}).
\end{align*}
\end{thm}
\begin{proof}
See Appendix.
\end{proof}

\begin{nas}\label{nas_h1_h2_e_n_c}
The spectral characteristics $\vec{h}_{\tau}^1(\lambda)$ and $\vec{h}_{\tau}^2(\lambda)$ of the optimal estimates $\widehat{B}\vec\zeta$ and $\widehat{A}\vec\eta$ of the functionals $B\vec\zeta$ and $A\vec\eta$ based on observations $\xi(t)+\eta(t)$ at points   $t<0$ are calculated by the formulas
\begin{align*}
\notag(\vec{h}_{\tau}^1(\lambda))^{\top}&=(\vec{B}_{\tau}(e^{i\lambda}))^{\top}
\frac{(1-e^{-i\lambda\tau})^d}{(i\lambda)^d}
\\
&\quad-
\frac{(-i\lambda)^d}{(1-e^{i\lambda \tau})^d}
\left(
\sum_{j=0}^{\infty}(\me P_{\tau}^{-1}D^{\tau}\me a)_j e^{i\lambda j}
\right)^{\top}(f(\lambda)+{\lambda}^{2d}g(\lambda)^{-1},
\\
\notag(\vec{h}_{\tau}^2(\lambda))^{\top}&=
(\vec{A}_{\tau}(e^{i\lambda }))^{\top}
 {g(\lambda)}\frac{(-i\lambda)^d}{(1-e^{i\lambda \tau})^d}(f(\lambda)+{\lambda}^{2d}g(\lambda))^{-1}
 \\
 &\quad-
\frac{(-i\lambda)^d}{(1-e^{i\lambda \tau})^d}
\left(
\sum_{j=0}^{\infty}(\me P_{\tau}^{-1}\me T_{\tau}\me a^{\tau})_j e^{i\lambda j}
\right)^{\top}(f(\lambda)+{\lambda}^{2d}g(\lambda))^{-1},
\end{align*}
respectively.
\end{nas}

\par

The optimal estimate
$\widehat{A}_{NT}\xi$ of the functional $A_{NT}\xi$ which depend on the unknown values of the process
$\xi(t)$  at points  $t\in[0,T(N+1)]$, based on observations of the process
$\xi(t)+\eta(t)$ at  points $t<0$ can be obtained by using Theorem \ref{thm1_e_n_c}.

\noindent We first formulate the following corollaries from Lemma \ref{predst A_cont} and Lemma \ref{predst H_cont}.

\begin{nas}\label{nas predst A_T_e_c}
The linear functional
\[A_{NT}\xi=\int_{0}^{(N+1)T}a(t)\xi(t)dt\]
allows the representation
\[
    A_{NT}\xi=B_{NT}\xi-V_{NT}\xi,\]
where
\[
    B_{NT}\xi=\int_0^{(N+1)T} b^{\tau,N}(t)\xi^{(d)}(t,\tau T)dt,\quad V_{NT}\xi=\int_{-\tau T d}^{0}v^{\tau,N}(t)\xi(t)dt,\]
and
\begin{align}\label{koefv_N_cont}
    v^{\tau,N}(t)&=\sum_{l=\left \lceil -\frac{t}{\tau T}\right \rceil}^{\min\left \{\left [\frac{{(N+1)T}-t}{\tau T}\right ],d\right \}}(-1)^l{d \choose l}b^{\tau, N}(t+l\tau T),\quad t\in[-\tau T d;0),
\\
\label{koef_N b_cont}
    b^{\tau,N}(t)&=\sum_{k=0}^{\left[\frac{{(N+1)T}-t}{\tau T}\right]}a( t +\tau T k)d(k)= D^{\tau T,N}\me a(t).
t\in[0;{(N+1)T}],
\end{align}
The linear transformation $D^{\tau T,N}$ acts on an arbitrary function $x(t)$, $t\in[0;{(N+1)T}]$, as follows
\[
    D^{\tau T,N}\me x(t)=\sum_{k=0}^{\left[\frac{{(N+1)T}-t}{\tau T}\right]}x(t+\tau T k)d(k).\]
\end{nas}
\begin{nas}\label{nas predst H_T_e_c}
The functional
$$H_{NT}\xi=B_{NT}\zeta-A_{NT}\eta$$
\[
 B_{NT}\zeta=\int_0^{T(N+1)} b^{\tau,N}(t)\zeta^{(d)}(t,\tau T)dt,\quad A_{NT}\eta=\int_{0}^{T(N+1)}a(t)\eta(t)dt,\]
can be represented in the form
\[
 H_{NT}\xi= \sum_{j=0}^N
{(\vec{b}^{\tau,N}_j)}^{\top}\vec{\zeta}^{(d)}_j- \sum_{j=0}^N
{(\vec{a}_j)}^{\top}\vec{\eta}_j
=B_{N}\vec \zeta - V_{N}\vec \eta=:H_N\vec \xi,
\]
where the vector
\[
\vec{b}^{\tau,N}_j=(b^{\tau,N}_{kj},k=1,2,\dots)^{\top}=
(b^{\tau,N}_{1j},b^{\tau,N}_{3j},b^{\tau,N}_{2j},\dots,b^{\tau,N}_{2k+1,j},b^{\tau,N}_{2k,j},\dots)^{\top},
\]
with the entries
\[
 b^{\tau,N}_{kj}=\langle b^{\tau,N}_j,\widetilde{e}_k\rangle =
\frac{1}{\sqrt{T}} \int_0^T b^{\tau,N}_j(v)e^{-2\pi
i\{(-1)^k\left[\frac{k}{2}\right]\}v/T}dv.
\]
The coefficients $\{\vec{a}_j, j=0,1,\dots,N\}$ and  $\{\vec{b}^{\tau,N}_j, j=0,1,\dots,N\}$ are related as
\be \label{a_b_N_relation}
\vec{b}_j^{\tau,N}=\sum_{m=j}^{N}\mt{diag}_{\infty}(d_{\tau}(m-j))\vec{a}_m=(D^{\tau}_N{\me a_N})_j,  \quad j=0,1,\dots,N.
\ee
where
$\me a_N=((\vec{a}_0)^{\top},(\vec{a}_1)^{\top}, \ldots,(\vec{a}_N)^{\top},0,\ldots)^{\top}$,
$D^{\tau}_N$ is a linear transformation  determined by a matrix with the infinite dimension matrix entries
 $D^{\tau}_N(k,j), k,j\geq0$ such that $D^{\tau}_N(k,j)=\mt{diag}_{\infty}(d_{\tau}(j-k))$ if $0\leq k\leq j \leq N$ and $D^{\tau}_N(k,j)=\mt{diag}_{\infty}(0)$ for $0\leq j<k$ or $j,k>N$.
\end{nas}

Put $a(t)=0$,
$t>T(N+1)$.  Define vector coefficients $\{\vec{a}_j^{\tau,N}:0\leq j\leq N+\tau d\}$ by the formula
\[
    \vec{a}_j^{\tau,N}
    =\sum_{l=\max\ld\{\ld\lceil\frac{j-N}{\tau}\rd\rceil,0\rd\}}^{\min\ld\{\ld[\frac{j}{\tau}\rd],d\rd\}}
    (-1)^l{d \choose l}\vec{a}(j-\tau l),\quad 0\leq j\leq N+\tau d,
\]
and a vector
\[
\me a^{\tau}_N=((\vec{a}_0^{\tau,N})^{\top},(\vec{a}_1^{\tau,N})^{\top},(\vec{a}_2^{\tau,N})^{\top},\ldots,
(\vec{a}_{N+\tau d}^{\tau,N})^{\top},0\ldots)^{\top}.
\]
The following theorem holds true.

\begin{thm}\label{thm2AN_e_n_c}
Consider two uncorrelated
processes: a stochastic
process $\xi(t)$, $t\in \mr R$ with periodically stationary  increments, which determines a
generated stationary  $d$th increment sequence
$\vec{\xi}^{(d)}_j$ with the spectral density matrix $f(\lambda)=\{f_{kn}(\lambda)\}_{k,n=1}^{\infty}$,
and   a periodically stationary stochastic
process $\eta(t)$, $t\in \mr R$,  which determines
   a
generated stationary  sequence
$\vec{\eta}_j$ with the spectral density matrix $g(\lambda)=\{g_{kn}(\lambda)\}_{k,n=1}^{\infty}$.
Let the coefficients $\vec {a}_j$, $ \vec{b}^{\tau}_j$, $j=0,1,\ldots,N$ generated by the function $a(t)$, $0\leq t\leq T(N+1)$, satisfy conditions
\begin{equation*} \label{coeff_a_N_e_n_c}
 \|\vec{a}_j\|<\infty, \quad \|\vec{a}_j\|^2=\sum_{k=1}^\infty |{a}_{kj}|^2,\quad    j=0,1,\dots,N,
\end{equation*}
and
\begin{equation*} \label{coeff_b_N_e_n_c}
 \|\vec{b}^{\tau,N}_j\|<\infty, \quad \|\vec{b}^{\tau,N}_j\|^2=\sum_{k=1}^\infty |{b}^{\tau,N}_{kj}|^2,\quad    j=0,1,\dots,N.
\end{equation*}
Let  minimality condition
(\ref{umova11_e_n_c}) be satisfied. The
optimal linear estimate $\widehat{A}_{NT}\xi$ of the functional
$A_{NT}\xi$  based on
observations of the process $\xi(t)+\eta(t)$ at points   $t<0$
is calculated by  formula
\begin{equation*} \label{otsinka A_N_e_n_c}
 \widehat{A}_{NT}\xi=\ip
(\vec{h}_{\tau,N}(\lambda))^{\top}d\vec{Z}_{\xi^{(d)}+\eta^{(d)}}(\lambda)-\int_{-\tau T d}^{0}v_{\tau,N}(t)\zeta(t)dt,
\end{equation*}
where the spectral characteristic $\vec{h}_{\tau,N}(\lambda)$ of the optimal
estimate $\widehat{A}_{NT}\xi$ is calculated by  formula
\begin{align*}
\notag(\vec{h}_{\tau,N}(\lambda))^{\top}
&=\frac{(1-e^{-i\lambda\tau})^d}{(i\lambda)^d}\left(\sum_{j=0}^{N}(D^{\tau}_N\me a_N)_je^{i\lambda j}\right)^{\top}
\\
\notag&\quad-
\frac{(-i\lambda)^d}{(1-e^{i\lambda \tau})^d}\left(\sum_{j=0}^{N+\tau d}\vec{a}^{\tau,N}_je^{i\lambda j}\right)^{\top}g(\lambda)
(f(\lambda)+{\lambda}^{2d}g(\lambda))^{-1}
\\
\notag&\quad-
\frac{(-i\lambda)^d}{(1-e^{i\lambda \tau})^d}
\left(
\sum_{j=0}^{\infty}(\me P_{\tau}^{-1}D_N^{\tau}\me a_N-\me P_{\tau}^{-1}\me T_{\tau}\me a_{\tau,N})_j e^{i\lambda j}
\right)^{\top}
\\
 &\quad\quad\times (f(\lambda)+{\lambda}^{2d}g(\lambda))^{-1},\label{spectr A_N_e_n_c}
\end{align*}
The value of the mean-square error
is calculated by  formula
\begin{align*}
\notag \Delta(f,g;\widehat{A}_{NT}\xi) &=\Delta(f,g;\widehat{H}_N\vec\xi)= \mt E|H_N\vec\xi-\widehat{H}_N\vec\xi|^2=
\\
\notag&=\langle D_N^{\tau}\me a_N- \me T_{\tau}\me
 a_{\tau,N},\me P_{\tau}^{-1}D_N^{\tau}\me a_N-\me P_{\tau}^{-1}\me T_{\tau}\me a_{\tau,N}\rangle
 \\
 &\quad+\langle\me Q_N\me a_N,\me  a_N\rangle, 
\end{align*}
 where
 $\me Q_{N}$ is a linear operator in the space $\ell_2$ defined by the matrix with the infinite dimensional  matrix  entries $(\me Q_{N})_{l,j}=Q_{l,j}$, $0\leq l,j\leq N$, and $(\me Q_{N})_{l,j}=0$ otherwise.
\end{thm}

\section{Prediction based on factorizations of spectral densities}\label{Prediction-factorizations}

Assume that the following canonical factorizations take place
\begin{equation} \label{fakt1_e_n_c}
 \frac{|1-e^{i\lambda\tau}|^{2d}}{\lambda^{2d}}
 (f(\lambda)+\lambda^{2d}g(\lambda))
 =\Theta_{\tau}(e^{-i\lambda})\Theta_{\tau}^*(e^{-i\lambda}),
 \\
 \Theta_{\tau}(e^{-i\lambda})=\sum_{k=0}^{\infty}\Theta_{\tau}(k)
  e^{-i\lambda k},
\end{equation}

\be \label{fakt3_e_n_c}
 g(\lambda)=\sum_{k=-\infty}^{\infty}g(k)e^{i\lambda k}=\Phi(e^{-i\lambda})\Phi^*(e^{-i\lambda}), \quad
 \Phi(e^{-i\lambda})=\sum_{k=0}^{\infty}\phi(k)e^{-i\lambda k}.
 \ee
Define the matrix-valued function
$\Psi_{\tau}(e^{-i\lambda})= \{\Psi_{ij}(e^{-i\lambda})\}_{i=\overline{1,M}}^{j=\overline{1,\infty}}$
by the equation
\[
\Psi_{\tau}(e^{-i\lambda})\Theta_{\tau}(e^{-i\lambda})=E_M,
\]
where $E_M$ is an identity $M\times M$ matrix.
Then the following factorization takes place
\begin{equation}
\frac{\lambda^{2d}}{|1-e^{i\lambda\tau}|^{2d}}
 (f(\lambda)+\lambda^{2d}g(\lambda))^{-1} =
 \Psi_{\tau}^*(e^{-i\lambda})\Psi_{\tau}(e^{-i\lambda}),
\,
 \Psi_{\tau}(e^{-i\lambda})=\sum_{k=0}^{\infty}\psi_{\tau}(k)e^{-i\lambda k},\label{fakt2_e_n_c}
\end{equation}

\begin{zau}\label{remark_density_adjoint_e_n_c}
Any spectral density matrix $f(\lambda)$ is self-adjoint: $f(\lambda)=f^*(\lambda)$. Thus, $(f(\lambda))^{\top}=\overline{f(\lambda)}$. One can check that an inverse spectral density $f^{-1}(\lambda)$ is also self-adjoint $f^{-1}(\lambda)=(f^{-1}(\lambda))^*$ and $(f^{-1}(\lambda))^{\top}=\overline{f^{-1}(\lambda)}$.
\end{zau}

The following Lemmas provide factorizations of the operators $\me P_{\tau}$ and $\me T_{\tau}$, which contain  coefficients of factorizations (\ref{fakt1_e_n_c}) -- (\ref{fakt2_e_n_c}).

\begin{lema}\label{lema_fact_2_e_n_c}
Let factorization (\ref{fakt1_e_n_c})  takes place and let $M\times \infty$ matrix function $\Psi_{\tau}(e^{-i\lambda})$ satisfy equation $\Psi_{\tau}(e^{-i\lambda})\Theta_{\tau}(e^{-i\lambda})=E_M$. Define the linear operators
 $ \Psi_{\tau}$ and  $ \Theta_{\tau}$ in the space $\ell_2$ by the matrices with the matrix entries
 $( \Psi_{\tau})_{k,j}=\psi_{\tau}(k-j)$, $( \Theta_{\tau})_{k,j}=\Theta_{\tau}(k-j)$ for $0\leq j\leq k$, $(
\Psi_{\tau})_{k,j}=0$, $(
\Theta_{\tau})_{k,j}=0$ for $0\leq k<j$.
Then:
\\
a) the linear operator $\me P_{\tau}$  admits the factorization \[\me
P_{\tau}=(\Psi_{\tau})^{\top} \overline{\Psi}_{\tau};\]
\\
b) the inverse operator $(\me
P_{\tau})^{-1}$ admits the factorization
\[
 (\me
P_{\tau})^{-1}= \overline{\Theta}_{\tau}(\Theta_{\tau})^{\top}.\]
\end{lema}

\begin{lema}\label{lema_fact_1_e_n_c}
 Let  factorizations (\ref{fakt1_e_n_c}) and (\ref{fakt3_e_n_c}) take place. Then the operator $\me T_{\tau}$ admits the representation
\[
\me T_{\tau}=(\Psi_{\tau})^{\top} \me Z_{\tau},
\]
 where $\me Z_{\tau}$ is a linear operator in the space $\ell_2$ defined by a matrix with the entries
\begin{eqnarray*}
 (\me Z_{\tau})_{k,j}&=&\sum_{l=j}^{\infty}\overline{\psi}_{\tau}(l-j)\overline{g}(l-k),\quad k,j\geq0
 \\
 g(k)&=&\sum_{m=\max\{0,-k\}}^{\infty}\phi(m)\phi^*(k+m),\quad k\in\mr Z.
\end{eqnarray*}
\end{lema}

\begin{zau}\label{remark_operator_PTa_e_n_c}
Lemma $\ref{lema_fact_2_e_n_c}$ and Lemma $\ref{lema_fact_1_e_n_c}$ imply the factorization
\[
(\me P_{\tau})^{-1}\me T_{\tau} \me a^{\tau}= \overline{\Theta}_{\tau}(\Theta_{\tau})^{\top}(\Psi_{\tau})^{\top} \me Z_{\tau}\me a^{\tau}
=\overline{\Theta}_{\tau} \me Z_{\tau}\me a^{\tau}=\overline{\Theta}_{\tau} \me e_{\tau},
\]
where $\me e_{\tau}:=\me Z_{\tau}\me a^{\tau}$.
\end{zau}

The proofs of Lemma \ref{lema_fact_2_e_n_c} and Lemma \ref{lema_fact_1_e_n_c}, as well as the justification of the   following representations of the spectral characteristics $\vec{h}_{\tau}^1(\lambda)$ and $\vec{h}_{\tau}^2(\lambda)$, correspond to the ones in \cite{Luz_Mokl_extra_noise_PCI} for the finite-dimensional vector stationary increment sequences.

The spectral characteristic $\vec{h}_{\tau}^2(\lambda)$ of the optimal estimate $\widehat{{A}}\vec \eta$ from Corollary \ref{nas_h1_h2_e_n_c} can be presented as
\[
 \vec h_{\tau}^2(\lambda)=\frac{(1-e^{-i\lambda\tau})^d}{(i\lambda)^d}\Psi^{\top}_{\tau}(e^{-i\lambda })\vec C_{\tau,g} (e^{-i\lambda}),
\]
where $\overline{\psi}_{\tau}=(\overline{\psi}_{\tau}(0),\overline{\psi}_{\tau}(1),\overline{\psi}_{\tau}(2),\ldots)$,
\[\vec C_{\tau,g} (e^{-i\lambda})=\sum_{m=1}^{\infty}(\overline{\psi}_{\tau} \me C_{\tau,g} )_m e^{-i\lambda m},\]
\begin{align*}
(\overline{\psi}_{\tau} \me C_{\tau,g} )_m&=\sum_{k=0}^{\infty}\overline{\psi}_{\tau}(k)\me c_{\tau,g}(k+m),
\\
\me c_{\tau,g}(m)&=\sum_{k=0}^{\infty}\overline{g}(m+k)\vec a^{\tau}_k
=\sum_{l=0}^{\infty}\overline{\phi}(l)\sum_{k=0}^{\infty}\phi^{\top}(l+m+k)\vec a^{\tau}_k
\\
&=\sum_{l=0}^{\infty}\overline{\phi}(l)(\widetilde{\Phi}\me a^{\tau})_{l+m},
\\
(\widetilde{\Phi}\me a^{\tau})_m&=\sum_{k=0}^{\infty}\phi^{\top}(m+k)\vec a^{\tau}_k.
\end{align*}

The spectral characteristic $\vec{h}_{\tau}^2(\lambda)$ of the optimal estimate $\widehat{ B}\vec\xi$ from Corollary \ref{nas_h1_h2_e_n_c} can be presented as
\begin{align*}
\vec{h}_{\tau}^1(\lambda)&=\frac{(1-e^{-i\lambda\tau})^d}{(i\lambda)^d}
  \ld(\vec B_{\tau}(e^{i\lambda})
 -\Psi^{\top}_{\tau}(e^{-i\lambda })\vec r_{\tau} (e^{i\lambda})\rd)
   \\&=\frac{(1-e^{-i\lambda\tau})^d}{(i\lambda)^d}\Psi^{\top}_{\tau}(e^{-i\lambda })\vec C_{\tau,1} (e^{-i\lambda}),
\end{align*}
where
\begin{align*}
\vec r_{\tau} (e^{i\lambda})&=
\sum_{m=0}^{\infty}( \theta^{\top}_{\tau}D^{\tau}\me A)_m e^{i\lambda m}
 =
\sum_{m=0}^{\infty}\sum_{p=0}^{\infty}\theta^{\top}_{\tau}(p) \vec b^{\tau}_{p+m} e^{i\lambda m},
\\
\vec C_{\tau,1} (e^{-i\lambda})&=
\sum_{m=1}^{\infty}( \theta^{\top}_{\tau}\widetilde{\me B}_{\tau} )_m e^{-i\lambda m}
 =
\sum_{m=1}^{\infty}\sum_{p=m}^{\infty}\theta^{\top}_{\tau}(p)\vec b^{\tau}_{p-m} e^{-i\lambda m}
 \\ &=
\sum_{m=1}^{\infty}\sum_{p=0}^{\infty}\theta^{\top}_{\tau}(m+p)\vec b^{\tau}_p e^{-i\lambda m},
\end{align*}
a vector $\theta^{\top}_{\tau}=((\Theta_{\tau}(0))^{\top},(\Theta_{\tau}(1))^{\top},(\Theta_{\tau}(2))^{\top},\ldots)$; $\me A$ is a linear symmetric operator determined by the matrix with the vector entries $(\me A)_{k,j}=\vec a_{k+j}$, $k,j\geq0$; $\widetilde{\me B}_{\tau}$ is a linear operator, which is determined by a matrix with the vector entries
 $( \widetilde{\me B}_{\tau})_{k,j}=\vec b^{\tau}_{k-j}$ for $0\leq j\leq k$, $(
\widetilde{\me B}_{\tau})_{k,j}=0$ for $0\leq k<j$.

Then the spectral characteristic $\vec h_{\tau}(\lambda)$ of the estimate $\widehat{A}\xi$ can be calculated by the formula
\begin{align}
 \notag \vec h_{\tau}(\lambda)&=\frac{(1-e^{-i\lambda\tau})^d}{(i\lambda)^d}
\ld(\sum_{k=0}^{\infty}\psi^{\top}_{\tau}(k) e^{-i\lambda k}\rd)
\sum_{m=1}^{\infty}\ld( \theta^{\top}_{\tau}\widetilde{\me B}_{\tau} -\overline{\psi}_{\tau} \me C_{\tau,g}\rd)_m e^{-i\lambda m}
\\
\notag&=\frac{(1-e^{-i\lambda\tau})^d}{(i\lambda)^d}\Psi^{\top}_{\tau}(e^{-i\lambda })\ld(\vec C_{\tau,1} (e^{-i\lambda})-\vec C_{\tau,g} (e^{-i\lambda})\rd)
 \\&=\vec B_{\tau}(e^{i\lambda})\frac{(1-e^{-i\lambda\tau})^d}{(i\lambda)^d}-\widetilde{h}_{\tau}(\lambda),
 \label{simple_spectr A_e_n_c}\end{align}
where
\begin{align*}
\widetilde{\vec h}_{\tau}(\lambda)&=\frac{(1-e^{-i\lambda\tau})^d}{(i\lambda)^d}
\ld(\sum_{k=0}^{\infty}\psi^{\top}_{\tau}(k) e^{-i\lambda k}\rd)\\
& \quad\times\ld(\sum_{m=0}^{\infty}( \theta^{\top}_{\tau}D^{\tau}\me A)_m e^{i\lambda m}+\sum_{m=1}^{\infty}(\overline{\psi}_{\tau} \me C_{\tau,g} )_m e^{-i\lambda m}\rd)
\\&=
\frac{(1-e^{-i\lambda\tau})^d}{(i\lambda)^d}
\Psi^{\top}_{\tau}(e^{-i\lambda })\ld(\vec r_{\tau} (e^{i\lambda})+\vec C_{\tau,g} (e^{-i\lambda})\rd),
\end{align*}

The value of the mean square error of the estimate $\widehat{A}\xi$ is calculated by the formula
\begin{align}
 \notag \Delta(f,g;\widehat{A}\xi)&=\Delta(f,g;\widehat{H}\vec\xi)= \mt E|H\vec\xi-\widehat{H}\vec\xi|^2
 \\\notag&=\frac{1}{2\pi}\int_{-\pi}^{\pi}(\vec A(e^{i\lambda}))^{\top}g(\lambda)\overline{\vec A(e^{i\lambda})}d\lambda
\\ \notag &\quad+
 \frac{1}{2\pi}\int_{-\pi}^{\pi}(\widetilde{\vec h}_{\tau}(e^{i\lambda}))^{\top}(f(\lambda)
 +\lambda^{2d}g(\lambda))\overline{\widetilde{\vec h}_{\tau}(e^{i\lambda})}d\lambda
\\ \notag &\quad-\frac{1}{2\pi}\int_{-\pi}^{\pi}\frac{(i\lambda)^d}{(1-e^{-i\lambda\tau})^d}(\widetilde{h}_{\tau}(e^{i\lambda}))^{\top}
g(\lambda)\overline{\vec A^{\tau}(e^{i\lambda})}d\lambda
\\\notag &\quad -
 \frac{1}{2\pi}\int_{-\pi}^{\pi}\frac{(-i\lambda)^d}{(1-e^{i\lambda \tau})^d}(\vec A^{\tau}(e^{i\lambda}))^{\top}
g(\lambda)
\overline{\widetilde{\vec h}_{\tau}(e^{i\lambda})}d\lambda
  \\\notag
 &=\|\Phi^{\top}\me a^{\tau}\|^2+\|\widetilde{\Phi}\me a^{\tau}\|_1^2
 +\ld\langle \theta^{\top}_{\tau}D^{\tau}\me A-\overline{\psi}_{\tau} \me C_{\tau,g},\theta^{\top}_{\tau}D^{\tau}\me A\rd\rangle
 \\ &\quad -\ld\langle\theta^{\top}_{\tau}D^{\tau}\me A,\me Z_{\tau}\me a^{\tau}\rd\rangle
 -\ld\langle \me Z_{\tau}\me a^{\tau},\overline{\psi}_{\tau} \me C_{\tau,g}\rd\rangle_1,
\label{simple_poh A_e_n_c}
\end{align}
where $\|\vec x\|_1^2=\langle \vec x, \vec x\rangle_1$, $\langle \vec x, \vec y\rangle_1=\sum_{j=1}^{\infty}(\vec x_j)^{\top}\overline{\vec y}_j$ for the vectors
\[\vec x=((\vec x_0)^{\top},(\vec x_1)^{\top},(\vec x_2)^{\top},\ldots)^{\top}, \quad \vec y=((\vec y_0)^{\top},(\vec y_1)^{\top},(\vec y_2)^{\top},\ldots)^{\top}.\]

The obtained results are summarized in the form of the following theorem.

\begin{thm}\label{thm3_e_n_c}
Let the conditions of Theorem \ref{thm1_e_n_c} be fulfilled  and the spectral densities $f(\lambda)$ and $g(\lambda)$ of the stochastic processes $ \xi(t)$ and $\eta(t)$ admit  canonical factorizations (\ref{fakt1_e_n_c}) -- (\ref{fakt2_e_n_c}). Then the spectral characteristic
$\vec h_{\tau}(\lambda)$ and the value of the mean square error $\Delta(f,g;\widehat{A}\xi)$ of the optimal estimate $\widehat{A}\xi$ of the functional $A\xi$ based on observations of the processes $ \xi(t)+\eta(t)$ at points $t<0$ can be calculated by formulas (\ref{simple_spectr A_e_n_c}) and (\ref{simple_poh A_e_n_c}) respectively.
\end{thm}

\section{Minimax (robust) method of prediction}\label{minimax_extrapolation}

Consider the estimation problem for the functional ${A}\xi$ based on the observations
$\xi(t)+\eta(t)$ at points $ t<0$ when the spectral densities of sequences  are not exactly known while a set $\md D=\md D_f\times\md D_g$ of admissible spectral densities is defined.
The minimax (robust) approach of
estimation is applied. It is formalized by the following two definitions.

\begin{ozn}
For a given class of spectral densities $\mathcal{D}=\md
D_f\times\md D_g$ the spectral densities
$f^0(\lambda)\in\mathcal{D}_f$, $g^0(\lambda)\in\md D_g$ are called
least favorable in the class $\mathcal{D}$ for the optimal linear
prediction of the functional $A\xi$  if the following relation holds
true:
\[\Delta(f^0,g^0)=\Delta(h(f^0,g^0);f^0,g^0)=
\max_{(f,g)\in\mathcal{D}_f\times\md
D_g}\Delta(h(f,g);f,g).\]
\end{ozn}

\begin{ozn}
For a given class of spectral densities $\mathcal{D}=\md
D_f\times\md D_g$ the spectral characteristic $h^0(\lambda)$ of
the optimal linear estimate of the functional $A\xi$ is called
minimax-robust if there are satisfied the conditions
\[h^0(\lambda)\in H_{\mathcal{D}}=\bigcap_{(f,g)\in\mathcal{D}_f\times\md D_g}L_2^{0-}(f(\lambda)+\lambda^{2d}g(\lambda)),\]
and
\[\min_{h\in H_{\mathcal{D}}}\max_{(f,g)\in \mathcal{D}_f\times\md D_g}\Delta(h;f,g)=\max_{(f,g)\in\mathcal{D}_f\times\md
D_g}\Delta(h^0;f,g).\]
\end{ozn}

Taking into account the introduced definitions and the derived relations we can verify that the following lemmas hold true.

\begin{lema}
Spectral densities $f^0\in\mathcal{D}_f$,
$g^0\in\mathcal{D}_g$ which satisfy condition (\ref{umova11_e_n_c})
are least favorable in the class $\md D=\md D_f\times\md D_g$ for
the optimal linear prediction of the functional $A\vec\xi$ if operators
$\me P_{\tau}^0$, $\me T_{\tau}^0$, $\me Q^0$  defined by the
Fourier coefficients of the functions
\[
\frac{\lambda^{2d}}{|1-e^{i\lambda\tau}|^{2d}}
\ld[g^0(\lambda)(f^0(\lambda)+{\lambda}^{2d}g^0(\lambda))^{-1}\rd]^{\top},
\]
\[
\dfrac{\lambda^{2d}}{|1-e^{i\lambda\tau}|^{2d}}\ld[(f^0(\lambda)+{\lambda}^{2d}g^0(\lambda))^{-1}\rd]^{\top},\]
\[
\ld[f^0(\lambda)(f^0(\lambda)+{\lambda}^{2d}g^0(\lambda))^{-1}g^0(\lambda)\rd]^{\top}
\]
determine a solution of the constrained optimization problem
\begin{multline}
\max_{(f,g)\in \mathcal{D}_f\times\md D_g}(\langle D^{\tau}\me a- \me T_{\tau}\me
 a^{\tau},\me P_{\tau}^{-1}D^{\tau}\me a-\me P_{\tau}^{-1}\me T_{\tau}\me a^{\tau}\rangle+\langle\me Q\me a,\me a\rangle)
\\
    = \langle D^{\tau}\me a- \me T^0_{\tau}\me
    a^{\tau},(\me P^0_{\tau})^{-1}D^{\tau}\me a-(\me P^0_{\tau})^{-1}\me T^0_{\tau}\me a^{\tau}\rangle+\langle\me Q^0\me a,\me a\rangle.
\label{minimax1_e_n_c}
\end{multline}
The minimax spectral characteristic $h^0=h_{\tau}(f^0,g^0)$ is calculated by formula (\ref{spectr_A_e_n_c}) if
$h_{\tau}(f^0,g^0)\in H_{\mathcal{D}}$.
\end{lema}

\begin{lema} The spectral densities $f^0\in\mathcal{D}_f$,
$g^0\in\mathcal{D}_g$ which admit  canonical factorizations (\ref{dd_cont}), (\ref{fakt1_e_n_c}) and (\ref{fakt3_e_n_c})
are least favourable densities in the class $\mathcal{D}$ for the optimal linear prediction
of the functional $A\xi$ based on observations of the process $\xi(t)+\eta(t)$ at points $t<0$ if the matrix coefficients
of  canonical factorizations
(\ref{fakt1_e_n_c}) and (\ref{fakt3_e_n_c})
determine a solution to the constrained optimization problem
\begin{align}
\notag
\|\Phi^{\top}\me a^{\tau}\|^2 & +\|\widetilde{\Phi}\me a^{\tau}\|_1^2
 +\ld\langle \theta^{\top}_{\tau}D^{\tau}\me A-\overline{\psi}_{\tau} \me C_{\tau,g},\theta^{\top}_{\tau}D^{\tau}\me A\rd\rangle
 \\&-  \ld\langle\theta^{\top}_{\tau}D^{\tau}\me A,\me Z_{\tau}\me a^{\tau}\rd\rangle
 -\ld\langle \me Z_{\tau}\me a^{\tau},\overline{\psi}_{\tau} \me C_{\tau,g}\rd\rangle_1\rightarrow\sup, \label{simple_minimax1_e_n_c}
\end{align}
for
\begin{align*}
 f(\lambda)&=\frac{\lambda^{2d}}{|1-e^{i\lambda\tau}|^{2d}}
\Theta_{\tau}(e^{-i\lambda})\Theta_{\tau}^*(e^{-i\lambda})
-\lambda^{2d}\Phi(e^{-i\lambda})\Phi^*(e^{-i\lambda})\in \mathcal{D}_f,
\\
 g(\lambda)&=\Phi(e^{-i\lambda})\Phi^*(e^{-i\lambda})\in \mathcal{D}_g.
\end{align*}
The minimax spectral characteristic $\vec h^0=\vec h_{\tau}(f^0,g^0)$ is calculated by formula (\ref{simple_spectr A_e_n_c}) if
$\vec h_{\tau}(f^0,g^0)\in H_{\mathcal{D}}$.
\end{lema}

\begin{lema} The spectral density $g^0\in\mathcal{D}_g$ which admits  canonical factorizations (\ref{fakt1_e_n_c}), (\ref{fakt3_e_n_c}) with the known spectral density $f(\lambda)$ is the least favourable in the class $\mathcal{D}_g$ for the optimal linear prediction
of the functional $A\xi$ based on observations of the process $\xi(t)+\eta(t)$ at points $t<0$ if the matrix coefficients
of the canonical factorizations
\begin{align} \label{fakt24_1_lf_e_n_c}
 f(\lambda)+\lambda^{2d}g^0(\lambda)&=\frac{\lambda^{2d}}{|1-e^{i\lambda\tau}|^{2d}}
\ld(\sum_{k=0}^{\infty}\theta^0_{\tau}(k)e^{-i\lambda k}\rd)\ld(\sum_{k=0}^{\infty}\theta^0_{\tau}(k)e^{-i\lambda k}\rd)^*,
\\ \label{fakt24_2_lf_e_n_c}
 g^0(\lambda)&=\ld(\sum_{k=0}^{\infty}\phi^0(k)e^{-i\lambda k}\rd)\ld(\sum_{k=0}^{\infty}\phi^0(k)e^{-i\lambda k}\rd)^*
\end{align}
and the equation $\Psi^0_{\tau}(e^{-i\lambda})\Theta^0_{\tau}(e^{-i\lambda})=E_M$ determine a solution to the constrained optimization problem
\begin{align}
\notag
\|\Phi^{\top}\me a^{\tau}\|^2&+\|\widetilde{\Phi}\me a^{\tau}\|_1^2
 +\ld\langle \theta^{\top}_{\tau}D^{\tau}\me A-\overline{\psi}_{\tau} \me C_{\tau,g},\theta^{\top}_{\tau}D^{\tau}\me A\rd\rangle
 \\&-  \ld\langle\theta^{\top}_{\tau}D^{\tau}\me A,\me Z_{\tau}\me a^{\tau}\rd\rangle
 -\ld\langle \me Z_{\tau}\me a^{\tau},\overline{\psi}_{\tau} \me C_{\tau,g}\rd\rangle_1\rightarrow\sup,
 \label{simple_minimax2_e_n_c}
 \end{align}
 for
\[
 g(\lambda)=\Phi(e^{-i\lambda})\Phi^*(e^{-i\lambda})\in \mathcal{D}_g.\]
The minimax spectral characteristic $\vec h^0=\vec h_{\tau}(f,g^0)$ is calculated by formula (\ref{simple_spectr A_e_n_c}) if
$\vec h_{\tau}(f,g^0)\in H_{\mathcal{D}}$.
\end{lema}

\begin{lema} The spectral density $f^0\in\mathcal{D}_f$ which admits  canonical factorizations
(\ref{dd_cont}), (\ref{fakt1_e_n_c}) with the known spectral density $g(\lambda)$ is the least favourable spectral density in the class
 $\md D_f$ for the optimal linear prediction
of the functional $A\xi$ based on observations of the process $\xi(t)+\eta(t)$ at points $t<0$ if matrix coefficients
of the canonical factorization
\be \label{fakt2_lf_e_n_c}
f^0(\lambda)+\lambda^{2d}g(\lambda)=\frac{\lambda^{2d}}{|1-e^{i\lambda\tau}|^{2d}}
\ld(\sum_{k=0}^{\infty}\theta^0_{\tau}(k)e^{-i\lambda k}\rd)\ld(\sum_{k=0}^{\infty}\theta^0_{\tau}(k)e^{-i\lambda k}\rd)^*,\ee
and the equation $\Psi^0_{\tau}(e^{-i\lambda})\Theta^0_{\tau}(e^{-i\lambda})=E_M$ determine a solution to the constrained optimization problem
\begin{align}
\ld\langle \theta^{\top}_{\tau}D^{\tau}\me A-\overline{\psi}_{\tau} \me C_{\tau,g},\theta^{\top}_{\tau}D^{\tau}\me A\rd\rangle
- \ld\langle\theta^{\top}_{\tau}D^{\tau}\me A,\me Z_{\tau}\me a^{\tau}\rd\rangle
 -\ld\langle \me Z_{\tau}\me a^{\tau},\overline{\psi}_{\tau} \me C_{\tau,g}\rd\rangle_1\rightarrow\sup,
 \label{simple_minimax3_e_n_c}
 \end{align}
for
\begin{align*}
 f(\lambda)&=\frac{\lambda^{2d}}{|1-e^{i\lambda\tau}|^{2d}}
\Theta_{\tau}(e^{-i\lambda})\Theta_{\tau}^*(e^{-i\lambda})
-\lambda^{2d}\Phi(e^{-i\lambda})\Phi^*(e^{-i\lambda})\in \mathcal{D}_f
\end{align*}
for the fixed matrix coefficients $\{\phi(k):k\geq0\}$. The minimax spectral characteristic $\vec h^0=\vec h_{\tau}(f^0,g)$ is calculated by formula (\ref{simple_spectr A_e_n_c}) if
$\vec h_{\tau}(f^0,g)\in H_{\mathcal{D}}$.
\end{lema}

For more detailed analysis of properties of the least favorable spectral densities and minimax-robust spectral characteristics we observe that
the minimax spectral characteristic $h^0$ and the least favourable spectral densities $(f^0,g^0)$
form a saddle
point of the function $\Delta(h;f,g)$ on the set
$H_{\mathcal{D}}\times\mathcal{D}$.

The saddle point inequalities
\[\Delta(h;f^0,g^0)\geq\Delta(h^0;f^0,g^0)\geq\Delta(h^0;f,g)
\quad\forall f\in \mathcal{D}_f,\forall g\in \mathcal{D}_g,\forall
h\in H_{\mathcal{D}}\] hold true if $h^0=\vec h_{\tau}(f^0,g^0)$ and
$\vec h_{\tau}(f^0,g^0)\in H_{\mathcal{D}}$, where $(f^0,g^0)$  is a
solution of the  constrained optimisation problem
\be  \label{cond-extr1_e_n_c}
\widetilde{\Delta}(f,g)=-\Delta(\vec h_{\tau}(f^0,g^0);f,g)\to
\inf,\quad (f,g)\in \mathcal{D},\ee
where the functional $\Delta(\vec h_{\tau}(f^0,g^0);f,g)$ is calculated by the formula
\begin{align*}
&\Delta(\vec h_{\tau}(f^0,g^0);f,g)=
\\
&=\frac{1}{2\pi}\int_{-\pi}^{\pi}
\frac{\lambda^{2d}}{|1-e^{i\lambda\tau}|^{2d}}
(\me C^{f0}_{\tau}(e^{i\lambda}))^{\top}(f^0(\lambda)+{\lambda}^{2d}g^0(\lambda))^{-1}f(\lambda)
\\
&\quad\quad\times
(f^0(\lambda)+{\lambda}^{2d}g^0(\lambda))^{-1}
\overline{\me C^{f0}_{\tau}(e^{i\lambda})}
d\lambda
\\
&\quad+\frac{1}{2\pi}\int_{-\pi}^{\pi}
\frac{\lambda^{4n}}{|1-e^{i\lambda\tau}|^{4n}}
(\me C^{g0}_{\tau}(e^{i\lambda}))^{\top}(f^0(\lambda)+{\lambda}^{2d}g^0(\lambda))^{-1} g(\lambda)
\\
&\quad\quad\times
(f^0(\lambda)+{\lambda}^{2d}g^0(\lambda))^{-1}
\overline{\me C^{g0}_{\tau}(e^{i\lambda})}
d\lambda,
\end{align*}
where
\begin{align*}
{\me C^{f0}_{\tau}(e^{i\lambda})}
&:=
\overline{g^0(\lambda)}\vec{A}_{\tau}(e^{i\lambda}) +
\sum_{j=0}^{\infty}((\me P^0_{\tau})^{-1}D^{\tau}\me a-(\me P^0_{\tau})^{-1}\me T^0_{\tau}\me a^{\tau})_j e^{i\lambda j},
\\
{\me C}^{g0}_{\tau}(e^{i \lambda})
&:=
{|1-e^{i\lambda\tau}|^{2d}}\lambda^{-2d}\overline{f^0(\lambda)}\vec{A}(e^{i\lambda})
\\
&\quad-{(1-e^{-i\lambda\tau})^d}
\sum_{j=0}^{\infty}((\me P^0_{\tau})^{-1}D^{\tau}\me a-(\me P^0_{\tau})^{-1}\me T^0_{\tau}\me a^{\tau})_j e^{i\lambda j}.
\end{align*}
or it is calculated by the formula
\begin{align*}
&\Delta\ld(\vec h_{\tau}(f^0,g^0);f,g\rd)=
\\
&=\frac{1}{2\pi}\ip\frac{|1-e^{i\lambda\tau}|^{2d}}{\lambda^{2d}}
(\me r^0_{\tau,f}(e^{-i\lambda}))^{\top}\Psi^0_{\tau}(e^{-i\lambda })f(\lambda)
(\Psi^0_{\tau}(e^{-i\lambda }))^*\overline{\me r^0_{\tau,f}(e^{-i\lambda})}d\lambda
\\
&\quad +
\frac{1}{2\pi}\ip
(\me r^0_{\tau,g}(e^{-i\lambda}))^{\top}\Psi^0_{\tau}(e^{-i\lambda })g(\lambda)
(\Psi^0_{\tau}(e^{-i\lambda }))^*\overline{\me r^0_{\tau,g}(e^{-i\lambda})}d\lambda,
\end{align*}
where
\begin{align*}
 \me r^0_{\tau,f}(e^{-i\lambda})&=\sum_{m=0}^{\infty}( (\theta^0_{\tau})^{\top}D^{\tau}\me A)_m e^{i\lambda m}+\sum_{m=1}^{\infty}(\overline{\psi}^0_{\tau} \me C^0_{\tau,g} )_m e^{-i\lambda m},
\\
 \me r^0_{\tau,g}(e^{-i\lambda})&=(1-e^{-i\lambda\tau})^{d}\ld(\sum_{m=0}^{\infty}( (\theta^0_{\tau})^{\top}D^{\tau}\me A)_m e^{i\lambda m}+\sum_{m=1}^{\infty}(\overline{\psi}_{\tau}^0 \me C^0_{\tau,g} )_m e^{-i\lambda m}\rd)
\\
&\quad-(\Theta^0_{\tau}(e^{-i\lambda}))^{\top}A(e^{i\lambda}).
\end{align*}

The constrained optimization problem (\ref{cond-extr1_e_n_c}) is equivalent to the unconstrained optimisation problem
\be  \label{uncond-extr}
\Delta_{\mathcal{D}}(f,g)=\widetilde{\Delta}(f,g)+ \delta(f,g|\mathcal{D}_f\times
\mathcal{D}_g)\to\inf,\ee
 where $\delta(f,g|\mathcal{D}_f\times
\mathcal{D}_g)$ is the indicator function of the set
$\mathcal{D}=\mathcal{D}_f\times\mathcal{D}_g$.
 Solution $(f^0,g^0)$ of this unconstrained optimization problem is characterized by the condition $0\in
\partial\Delta_{\mathcal{D}}(f^0,g^0)$, where
$\partial\Delta_{\mathcal{D}}(f^0,g^0)$ is the subdifferential of the functional $\Delta_{\mathcal{D}}(f,g)$ at point $(f^0,g^0)\in \mathcal{D}=\mathcal{D}_f\times\mathcal{D}_g$, that is the set of all continuous linear functionals $\Lambda$ on $L_1\times L_1$ which satisfy the inequality
$\Delta_{\mathcal{D}}(f,g)-\Delta_{\mathcal{D}}(f^0,g^0)\geq \Lambda ((f,g)-(f^0,g^0)), (f,g)\in \mathcal{D}$ (see  \cite{Moklyachuk2015,Rockafellar} for more details).
 This condition makes it possible to find the least favourable spectral densities in some special classes of spectral densities $\mathcal{D}=\mathcal{D}_f\times\mathcal{D}_g$.

The form of the functionals $\Delta(\vec h_{\tau}(f^0,g^0);f,g)$, $\Delta(\vec h_{\tau}(f^0,p^0);f,p)$ is convenient for application the Lagrange method of indefinite multipliers for
finding solution to the problem (\ref{uncond-extr}).
Making use of the method of Lagrange multipliers and the form of
subdifferentials of the indicator functions $\delta(f,g|\mathcal{D}_f\times
\mathcal{D}_g)$, $\delta(f,p|\mathcal{D}_f\times
\mathcal{D}_p)$ of the sets
$\mathcal{D}_f\times\mathcal{D}_g$, $\mathcal{D}_f\times\mathcal{D}_p$ of spectral densities,
we describe relations that determine least favourable spectral densities in some special classes
of spectral densities (see \cite{Luz_Mokl_book,Moklyachuk2015} for additional details).

\subsection{Least favorable spectral densities in classes $\md D_0 \times \md D_{1\delta}$}\label{set1_e_n_c}

Consider the prediction problem for the functional $A{\xi}$
 which depends on unobserved values of a process $\xi(t)$ with stationary increments based on
 observations of the process $\xi(t)+\eta(t)$ at points $t<0$ under the condition that the sets of admissible spectral densities
 $\md D_{0}^k, {D_{1\delta}^{k}},k=1,2,3,4$ are defined as follows:
\begin{align*}
\md D_{0}^{1} &= \bigg\{f(\lambda )\left|\frac{1}{2\pi} \int
_{-\pi}^{\pi}
\frac{|1-e^{i\lambda\tau}|^{2d}}{\lambda^{2d}}
f(\lambda )d\lambda  =P\right.\bigg\},
\\
 \md D_{0}^{2} &=\bigg\{f(\lambda )\left|\frac{1}{2\pi }
\int _{-\pi }^{\pi}
\frac{|1-e^{i\lambda\tau}|^{2d}}{\lambda^{2d}}
{\rm{Tr}}\,[ f(\lambda )]d\lambda =p\right.\bigg\},
\\
\md D_{0}^{3} &=\bigg\{f(\lambda )\left|\frac{1}{2\pi }
\int _{-\pi}^{\pi}
\frac{|1-e^{i\lambda\tau}|^{2d}}{\lambda^{2d}}
f_{kk} (\lambda )d\lambda =p_{k}, k=\overline{1,\infty}\right.\bigg\},
\\
\md D_{0}^{4} &=\bigg\{f(\lambda )\left|\frac{1}{2\pi} \int _{-\pi}^{\pi}
\frac{|1-e^{i\lambda\tau}|^{2d}}{\lambda^{2d}}
\left\langle B_{1} ,f(\lambda )\right\rangle d\lambda  =p\right.\bigg\},
\end{align*}
and
\begin{align*}
\md D_{1\delta}^{1}&=\left\{g(\lambda )\biggl|\frac{1}{2\pi} \int_{-\pi}^{\pi}
\left|g_{ij} (\lambda )-g_{ij}^{1} (\lambda)\right|d\lambda  \le \delta_{i}^j, i,j=\overline{1,\infty}\right\}.
\\
\md D_{1\delta}^{2}&=\left\{g(\lambda )\biggl|\frac{1}{2\pi} \int_{-\pi}^{\pi}
\left|{\rm{Tr}}(g(\lambda )-g_{1} (\lambda))\right|d\lambda \le \delta\right\};
\\
\md D_{1\delta}^{3}&=\left\{g(\lambda )\biggl|\frac{1}{2\pi } \int_{-\pi}^{\pi}
\left|g_{kk} (\lambda )-g_{kk}^{1} (\lambda)\right|d\lambda  \le \delta_{k}, k=\overline{1,\infty}\right\};
\\
\md D_{1\delta}^{4}&=\left\{g(\lambda )\biggl|\frac{1}{2\pi } \int_{-\pi}^{\pi}
\left|\left\langle B_{2} ,g(\lambda )-g_{1}(\lambda )\right\rangle \right|d\lambda  \le \delta\right\};
\end{align*}

\noindent
Here  $g_{1} ( \lambda )=\{g_{ij}^{1} ( \lambda )\}_{i,j=1}^{\infty}$ is a fixed spectral density, $p,  p_k, k=\overline{1,\infty}$, are given numbers, $P$, $B_2$ are a given positive-definite Hermitian matrices,
$\delta,\delta_{k},k=\overline{1,\infty}$, $\delta_{i}^{j}, i,j=\overline{1,\infty}$, are given numbers.

From the condition $0\in\partial\Delta_{\mathcal{D}}(f^0,g^0)$
we find the following equations which determine the least favourable spectral densities for these given sets of admissible spectral densities.

For the first set of admissible spectral densities $\md D_{0}^1 \times\md D_{1\delta}^{1}$,   we have equations
\begin{multline} \label{eq_4_1f_e_n_c}
\left(
{\me C^{f0}_{\tau}(e^{i\lambda})}
\right)
\left(
{\me C^{f0}_{\tau}(e^{i\lambda})}
\right)^{*}=
\\
=\left(\frac{|1-e^{i\lambda\tau}|^{2d}}{\lambda^{2d}} (f^0(\lambda)+{\lambda}^{2d}g^0(\lambda))\right)
\vec{\alpha}_f\cdot \vec{\alpha}_f^{*}
\times
\\
\times
\left(\frac{|1-e^{i\lambda\tau}|^{2d}}{\lambda^{2d}} (f^0(\lambda)+{\lambda}^{2d}g^0(\lambda))\right),
\end{multline}
\begin{multline}  \label{eq_5_1g_e_n_c}
\left(
{\me C}^{g0}_{\tau}(e^{i \lambda})
\right)
\left(
{\me C}^{g0}_{\tau}(e^{i \lambda})
\right)^{*}=
\\
=
\left(\frac{|1-e^{i\lambda\tau}|^{2d}}{\lambda^{2d}}(f^0(\lambda)+{\lambda}^{2d}g^0(\lambda))\right)
\left \{ \beta_{ij}\gamma_{ij} ( \lambda ) \right \}_{i,j=1}^{\infty}
\times
\\
\times
\left(\frac{|1-e^{i\lambda\tau}|^{2d}}{\lambda^{2d}}(f^0(\lambda)+{\lambda}^{2d}g^0(\lambda))\right),
\end{multline}
\begin{equation} \label{eq_5_1c_e_n_c}
\frac{1}{2 \pi} \int_{- \pi}^{ \pi} \left|g^0_{ij}(\lambda)-g_{ij}^{1}( \lambda ) \right|d\lambda = \delta_{i}^{j},
\end{equation}
where $\vec{\alpha}_f$, $ \beta_{ij}$ are Lagrange multipliers,  functions $\left| \gamma_{ij} ( \lambda ) \right| \le 1$ and
\[
\gamma_{ij} ( \lambda )= \frac{g_{ij}^{0} ( \lambda )-g_{ij}^{1} (\lambda )}{ \left|g_{ij}^{0} ( \lambda )-g_{ij}^{1}(\lambda) \right|}: \; g_{ij}^{0} ( \lambda )-g_{ij}^{1} ( \lambda ) \ne 0, \; i,j= \overline{1,\infty}.
\]

For the second set of admissible spectral densities $\md D_{0}^2 \times\md D_{1\delta}^{2}$,  we have equations
\be \label{eq_4_2f_e_n_c}
\left(
{\me C^{f0}_{\tau}(e^{i\lambda})}
\right)
\left(
{\me C^{f0}_{\tau}(e^{i\lambda})}
\right)^{*}=
\alpha_f^{2} \left(\frac{|1-e^{i\lambda\tau}|^{2d}}{\lambda^{2d}} (f^0(\lambda)+{\lambda}^{2d}g^0(\lambda))\right)^2,
\ee
\begin{equation} \label{eq_5_2g_e_n_c}
\left(
{\me C}^{g0}_{\tau}(e^{i \lambda})
\right)
\left(
{\me C}^{g0}_{\tau}(e^{i \lambda})
\right)^{*}
=
\beta^{2} \gamma_2( \lambda )\left(\frac{|1-e^{i\lambda\tau}|^{2d}}{\lambda^{2d}}(f^0(\lambda)+{\lambda}^{2d}g^0(\lambda))\right)^2,
\end{equation}
\begin{equation} \label{eq_5_2c_e_n_c}
\frac{1}{2 \pi} \int_{-\pi}^{ \pi}
\left|{\mathrm{Tr}}\, (g_0( \lambda )-g_{1}(\lambda )) \right|d\lambda =\delta,
\end{equation}

\noindent where $\alpha _{f}^{2}$, $ \beta^{2}$ are Lagrange multipliers,   the function $\left| \gamma_2( \lambda ) \right| \le 1$ and
\[\gamma_2( \lambda )={ \mathrm{sign}}\; ({\mathrm{Tr}}\, (g_{0} ( \lambda )-g_{1} ( \lambda ))): \; {\mathrm{Tr}}\, (g_{0} ( \lambda )-g_{1} ( \lambda )) \ne 0.\]

For the third set of admissible spectral densities $\md D_{0}^3 \times\md D_{1\delta}^{3}$,  we have equations
\begin{multline} \label{eq_4_3f_e_n_c}
\left(
{\me C^{f0}_{\tau}(e^{i\lambda})}
\right)
\left(
{\me C^{f0}_{\tau}(e^{i\lambda})}
\right)^{*}=
\\
=\left(\frac{|1-e^{i\lambda\tau}|^{2d}}{\lambda^{2d}} (f^0(\lambda)+{\lambda}^{2d}g^0(\lambda))\right)
\left\{\alpha _{fk}^{2} \delta _{kl} \right\}_{k,l=1}^{\infty}
\times
\\
\times
\left(\frac{|1-e^{i\lambda\tau}|^{2d}}{\lambda^{2d}} (f^0(\lambda)+{\lambda}^{2d}g^0(\lambda))\right),
\end{multline}
\begin{multline}   \label{eq_5_3g_e_n_c}
\left(
{\me C}^{g0}_{\tau}(e^{i \lambda})
\right)
\left(
{\me C}^{g0}_{\tau}(e^{i \lambda})
\right)^{*}=
\\
=
\left(\frac{|1-e^{i\lambda\tau}|^{2d}}{\lambda^{2d}}(f^0(\lambda)+{\lambda}^{2d}g^0(\lambda))\right)
\left \{ \beta_{k}^{2} \gamma^2_{k} ( \lambda ) \delta_{kl} \right \}_{k,l=1}^{\infty}
\times
\\
\times
\left(\frac{|1-e^{i\lambda\tau}|^{2d}}{\lambda^{2d}}(f^0(\lambda)+{\lambda}^{2d}g^0(\lambda))\right),
\end{multline}
\begin{equation} \label{eq_5_3c_e_n_c}
\frac{1}{2 \pi} \int_{- \pi}^{ \pi}  \left|g^0_{kk} ( \lambda)-g_{kk}^{1} ( \lambda ) \right| d\lambda =\delta_{k},
\end{equation}

\noindent where $\alpha _{fk}^{2}$, $\beta_{k}^{2}$ are Lagrange multipliers, $\delta _{kl}$ are Kronecker symbols,  functions $\left| \gamma^2_{k} ( \lambda ) \right| \le 1$ and
\[\gamma_{k}^2( \lambda )={ \mathrm{sign}}\;(g_{kk}^{0}( \lambda)-g_{kk}^{1} ( \lambda )): \; g_{kk}^{0} ( \lambda )-g_{kk}^{1}(\lambda ) \ne 0, \; k= \overline{1,\infty}.\]

For the fourth set of admissible spectral densities $\md D_{0}^4 \times\md D_{1\delta}^{4}$,  we have equations
\begin{multline} \label{eq_4_4f_e_n_c}
\left(
{\me C^{f0}_{\tau}(e^{i\lambda})}
\right)
\left(
{\me C^{f0}_{\tau}(e^{i\lambda})}
\right)^{*}=
\\
=
\alpha_f^{2} \left(\frac{|1-e^{i\lambda\tau}|^{2d}}{\lambda^{2d}} (f^0(\lambda)+{\lambda}^{2d}g^0(\lambda))\right)
B_{1}^{\top}
\times
\\
\times
\left(\frac{|1-e^{i\lambda\tau}|^{2d}}{\lambda^{2d}} (f^0(\lambda)+{\lambda}^{2d}g^0(\lambda))\right),
\end{multline}
\begin{multline}   \label{eq_5_4g_e_n_c}
\left(
{\me C}^{g0}_{\tau}(e^{i \lambda})
\right)
\left(
{\me C}^{g0}_{\tau}(e^{i \lambda})
\right)^{*}=
\\
=
\beta^{2} \gamma_2'( \lambda )
\left(\frac{|1-e^{i\lambda\tau}|^{2d}}{\lambda^{2d}}(f^0(\lambda)+{\lambda}^{2d}g^0(\lambda))\right)
B_{2}^{ \top}
\times
\\
\times
\left(\frac{|1-e^{i\lambda\tau}|^{2d}}{\lambda^{2d}}(f^0(\lambda)+{\lambda}^{2d}g^0(\lambda))\right),
\end{multline}
\begin{equation} \label{eq_5_4c_e_n_c}
\frac{1}{2 \pi} \int_{- \pi}^{ \pi}  \left| \left \langle B_{2}, g^0( \lambda )-g_{1} ( \lambda ) \right \rangle \right|d\lambda
= \delta,
\end{equation}
where $\alpha _{f}^{2}$,  $\beta^{2}$ are Lagrange multipliers,   function
$\left| \gamma_2' ( \lambda ) \right| \le 1$ and
\[\gamma_2' ( \lambda )={ \mathrm{sign}}\; \left \langle B_{2},g^{0} ( \lambda )-g_{1} ( \lambda ) \right \rangle : \; \left \langle B_{2},g^{0} ( \lambda )-g_{1} ( \lambda ) \right \rangle \ne 0.\]

The derived results are summarized in the  following theorem.

\begin{thm}
The least favorable spectral densities $f^{0}(\lambda), $ $g^{0}(\lambda), $ in the classes $\md D_{0}^k \times\md D_{1\delta}^{k},k=1,2,3,4$
for the optimal linear predictionof the functional  $A{\xi}$ from observations of the process ${\xi}(t)+{\eta}(t)$ at points  $t<0$  are determined by equations
\eqref{eq_4_1f_e_n_c}--\eqref{eq_5_1c_e_n_c},  \eqref{eq_4_2f_e_n_c}--\eqref{eq_5_2c_e_n_c}, \eqref{eq_4_3f_e_n_c}--\eqref{eq_5_3c_e_n_c}, \eqref{eq_4_4f_e_n_c}--\eqref{eq_5_4c_e_n_c},
respectively, the minimality condition (\ref{umova11_e_n_c}), the
  constrained optimization problem (\ref{minimax1_e_n_c}) and restrictions  on densities from the corresponding classes $\md D_{0}^k, \md D_{1\delta}^{k},k=1,2,3,4$.  The minimax-robust spectral characteristic of the optimal estimate of the functional $A{\xi}$ is determined by   formula (\ref{spectr_A_e_n_c}).
\end{thm}

In the case where the spectral densities $f(\lambda)$ and $g(\lambda)$ admit canonical factorizations  (\ref{fakt1_e_n_c}) and (\ref{fakt3_e_n_c}),   the    equation for the least favourable spectral densities are described below.

For the first set of admissible spectral densities $\md D_{0}^1 \times\md D_{1\delta}^{1}$:
\begin{align} \label{eq_4_1f_fact_e_n_c}
\left(
{\me r^{0}_{\tau,f}(e^{i\lambda})}
\right)
\left(
{\me r^{0}_{\tau,f}(e^{i\lambda})}
\right)^{*}
&=(\Theta_{\tau}(e^{-i\lambda}))^{\top}
\vec{\alpha}_f\cdot \vec{\alpha}_f^{*}
\overline{\Theta_{\tau}(e^{-i\lambda})},
\\  \label{eq_5_1g_fact_e_n_c}
\left(
{\me r^{0}_{\tau,g}(e^{i\lambda})}
\right)
\left(
{\me r^{0}_{\tau,g}(e^{i\lambda})}
\right)^{*}
&=
(\Theta_{\tau}(e^{-i\lambda}))^{\top}
\left \{ \beta_{ij}\gamma_{ij} ( \lambda ) \right \}_{i,j=1}^{T}
\overline{\Theta_{\tau}(e^{-i\lambda})},
\end{align}
\begin{equation} \label{eq_5_1c_fact_e_n_c}
\frac{1}{2 \pi} \int_{- \pi}^{ \pi} \left|g^0_{ij}(\lambda)-g_{ij}^{1}( \lambda ) \right|d\lambda = \delta_{i}^{j},
\end{equation}
where $\vec{\alpha}_f$, $ \beta_{ij}$ are Lagrange multipliers,  functions $\left| \gamma_{ij} ( \lambda ) \right| \le 1$ and
\[
\gamma_{ij} ( \lambda )= \frac{g_{ij}^{0} ( \lambda )-g_{ij}^{1} (\lambda )}{ \left|g_{ij}^{0} ( \lambda )-g_{ij}^{1}(\lambda) \right|}: \; g_{ij}^{0} ( \lambda )-g_{ij}^{1} ( \lambda ) \ne 0, \; i,j= \overline{1,\infty}.
\]

For the second set of admissible spectral densities $\md D_{0}^2 \times\md D_{1\delta}^{2}$:
\begin{align} \label{eq_4_2f_fact_e_n_c}
\left(
{\me r^{0}_{\tau,f}(e^{i\lambda})}
\right)
\left(
{\me r^{0}_{\tau,f}(e^{i\lambda})}
\right)^{*}&=
\alpha_f^{2} (\Theta_{\tau}(e^{-i\lambda}))^{\top}\overline{\Theta_{\tau}(e^{-i\lambda})},
\\ \label{eq_5_2g_fact_e_n_c}
\left(
{\me r^{0}_{\tau,g}(e^{i\lambda})}
\right)
\left(
{\me r^{0}_{\tau,g}(e^{i\lambda})}
\right)^{*}&=
\beta^{2} \gamma_2( \lambda )(\Theta_{\tau}(e^{-i\lambda}))^{\top}\overline{\Theta_{\tau}(e^{-i\lambda})},
\end{align}
\begin{equation} \label{eq_5_2c_fact_e_n_c}
\frac{1}{2 \pi} \int_{-\pi}^{ \pi}
\left|{\mathrm{Tr}}\, (g_0( \lambda )-g_{1}(\lambda )) \right|d\lambda =\delta,
\end{equation}

\noindent where $\alpha _{f}^{2}$, $ \beta^{2}$ are Lagrange multipliers,   function $\left| \gamma_2( \lambda ) \right| \le 1$ and
\[\gamma_2( \lambda )={ \mathrm{sign}}\; ({\mathrm{Tr}}\, (g_{0} ( \lambda )-g_{1} ( \lambda ))): \; {\mathrm{Tr}}\, (g_{0} ( \lambda )-g_{1} ( \lambda )) \ne 0.\]

For the third set of admissible spectral densities $\md D_{0}^3 \times\md D_{1\delta}^{3}$:
\begin{align} \label{eq_4_3f_fact_e_n_c}
\left(
{\me r^{0}_{\tau,f}(e^{i\lambda})}
\right)
\left(
{\me r^{0}_{\tau,f}(e^{i\lambda})}
\right)^{*}
&=(\Theta_{\tau}(e^{-i\lambda}))^{\top}
\left\{\alpha _{fk}^{2} \delta _{kl} \right\}_{k,l=1}^{T}
\overline{\Theta_{\tau}(e^{-i\lambda})},
\\   \label{eq_5_g_fact_e_n_c}
\left(
{\me r^{0}_{\tau,g}(e^{i\lambda})}
\right)
\left(
{\me r^{0}_{\tau,g}(e^{i\lambda})}
\right)^{*}
&=
(\Theta_{\tau}(e^{-i\lambda}))^{\top}
\left \{ \beta_{k}^{2} \gamma^2_{k} ( \lambda ) \delta_{kl} \right \}_{k,l=1}^{T}
\overline{\Theta_{\tau}(e^{-i\lambda})},
\end{align}
\begin{equation} \label{eq_5_3c_fact_e_n_c}
\frac{1}{2 \pi} \int_{- \pi}^{ \pi}  \left|g^0_{kk} ( \lambda)-g_{kk}^{1} ( \lambda ) \right| d\lambda =\delta_{k},
\end{equation}

\noindent where $\alpha_{fk}^{2}$, $\beta_{k}^{2}$ are Lagrange multipliers, $\delta _{kl}$ are Kronecker symbols,  functions $\left| \gamma^2_{k} ( \lambda ) \right| \le 1$ and
\[\gamma_{k}^2( \lambda )={ \mathrm{sign}}\;(g_{kk}^{0}( \lambda)-g_{kk}^{1} ( \lambda )): \; g_{kk}^{0} ( \lambda )-g_{kk}^{1}(\lambda ) \ne 0, \; k= \overline{1,\infty}.\]

For the fourth set of admissible spectral densities $\md D_{0}^4 \times\md D_{1\delta}^{4}$:

\begin{align} \label{eq_4_4f_fact_e_n_c}
\left(
{\me r^{0}_{\tau,f}(e^{i\lambda})}
\right)
\left(
{\me r^{0}_{\tau,f}(e^{i\lambda})}
\right)^{*}
&=
\alpha_f^{2} (\Theta_{\tau}(e^{-i\lambda}))^{\top}
B_{1}
\overline{\Theta_{\tau}(e^{-i\lambda})},
\\  \label{eq_5_3g_fact_e_n_c}
\left(
{\me r^{0}_{\tau,g}(e^{i\lambda})}
\right)
\left(
{\me r^{0}_{\tau,g}(e^{i\lambda})}
\right)^{*}
&=
\beta^{2} \gamma_2'( \lambda )
(\Theta_{\tau}(e^{-i\lambda}))^{\top}
B_{2}
\overline{\Theta_{\tau}(e^{-i\lambda})},
\end{align}
\begin{equation} \label{eq_5_4c_fact_e_n_c}
\frac{1}{2 \pi} \int_{- \pi}^{ \pi}  \left| \left \langle B_{2}, g^0( \lambda )-g_{1} ( \lambda ) \right \rangle \right|d\lambda
= \delta,
\end{equation}
where $\alpha_f^{2}$,  $\beta^{2}$ are Lagrange multipliers,   function
$\left| \gamma_2' ( \lambda ) \right| \le 1$ and
\[\gamma_2' ( \lambda )={ \mathrm{sign}}\; \left \langle B_{2},g^{0} ( \lambda )-g_{1} ( \lambda ) \right \rangle : \; \left \langle B_{2},g^{0} ( \lambda )-g_{1} ( \lambda ) \right \rangle \ne 0.\]

The following theorems  hold true.

\begin{thm}
 The least favorable spectral densities $f^{0}(\lambda)$,  $g^{0}(\lambda)$  in the classes $\md D_{0}^k \times\md D_{1\delta}^{k}$, $k=1,2,3,4$ for the optimal linear prediction of the functional  $A{\xi}$ from observations of the process $\xi(t)+\eta(t)$ at points  $t<0$  are determined
by  canonical factorizations   (\ref{fakt1_e_n_c}) and (\ref{fakt3_e_n_c}),
 equations
\eqref{eq_4_1f_fact_e_n_c}--\eqref{eq_5_1c_fact_e_n_c},  \eqref{eq_4_2f_fact_e_n_c}--\eqref{eq_5_2c_fact_e_n_c}, \eqref{eq_4_3f_fact_e_n_c}--\eqref{eq_5_3c_fact_e_n_c}, \eqref{eq_4_4f_fact_e_n_c}--\eqref{eq_5_4c_fact_e_n_c},
respectively,
 constrained optimization problem (\ref{simple_minimax1_e_n_c}) and restrictions  on densities from the corresponding classes $\md D_{0}^k, \md D_{1\delta}^{k},k=1,2,3,4$.
 The minimax-robust spectral characteristic of the optimal estimate of the functional $A\vec{\xi}$ is determined by  formula (\ref{simple_spectr A_e_n_c}).
\end{thm}

\begin{thm}
 In the case where the spectral density $g(\lambda)$ is known, the least favorable spectral density $f^{0}(\lambda)$ in the classes $\md D_{0}^k $, $k=1,2,3,4$
 for the optimal linear predictionof the functional  $A{\xi}$ from observations of the process $ {\xi}(t)+ {\eta}(t)$ at points  $t<0$  is determined
by  canonical factorizations (\ref{fakt3_e_n_c}) and (\ref{fakt1_e_n_c}),
 equations
\eqref{eq_4_1f_fact_e_n_c},  \eqref{eq_4_2f_fact_e_n_c}, \eqref{eq_4_3f_fact_e_n_c}, \eqref{eq_4_4f_fact_e_n_c},
respectively,
 constrained optimization problem (\ref{simple_minimax3_e_n_c}) and restrictions  on density from the corresponding classes $\md D_{0}^k$, $k=1,2,3,4$.  The minimax-robust spectral characteristic of the optimal estimate of the functional $A{\xi}$ is determined by   formula (\ref{simple_spectr A_e_n_c}).
\end{thm}

\subsection{Least favorable spectral densities in classes $\md D_{\varepsilon}\times \md  D_V^U$}\label{set2_e_n_c}

Consider the prediction problem for the functional $A{\xi}$
 depending on unobserved values of the process $\vec\xi(m)$ with stationary increments based on observations of the process $\xi(t)+\eta(t)$ at points $t<0$ under the condition that the sets of admissible spectral densities  $ \md D_{\varepsilon}^{k},  \md D_{V}^{Uk}, k=1,2,3,4$ are defined as follows:
\begin{equation*}
\md D_{\varepsilon }^{1}=\bigg\{f(\lambda )\bigg|f(\lambda)=(1-\varepsilon )f_{1} (\lambda )+\varepsilon W(\lambda ),
\frac{1}{2\pi } \int _{-\pi}^{\pi}
\frac{|1-e^{i\lambda\tau}|^{2d}}{\lambda^{2d}}
f(\lambda )d\lambda=P\bigg\},
\end{equation*}
\begin{multline*}
\md D_{\varepsilon }^{2}  =\bigg\{f(\lambda )\bigg|{\mathrm{Tr}}\,
[f(\lambda )]=(1-\varepsilon ) {\mathrm{Tr}}\,  [f_{1} (\lambda
)]+\varepsilon {\mathrm{Tr}}\,  [W(\lambda )],\\
\frac{1}{2\pi} \int _{-\pi}^{\pi}
\frac{|1-e^{i\lambda\tau}|^{2d}}{\lambda^{2d}}
{\mathrm{Tr}}\,
[f(\lambda )]d\lambda =p \bigg\};
\end{multline*}
\begin{multline*}
\md D_{\varepsilon }^{3}  =\bigg\{f(\lambda )\bigg|f_{kk} (\lambda)
=(1-\varepsilon )f_{kk}^{1} (\lambda )+\varepsilon w_{kk}(\lambda),\\
\frac{1}{2\pi} \int _{-\pi}^{\pi}
\frac{|1-e^{i\lambda\tau}|^{2d}}{\lambda^{2d}}
f_{kk} (\lambda)d\lambda  =p_{k} , k=\overline{1,\infty}\bigg\};
\end{multline*}
\begin{multline*}
\md D_{\varepsilon }^{4} =\bigg\{f(\lambda )\bigg|\left\langle B_{1},f(\lambda )\right\rangle =(1-\varepsilon )\left\langle B_{1},f_{1} (\lambda )\right\rangle+\varepsilon \left\langle B_{1},W(\lambda )\right\rangle,\\
\frac{1}{2\pi}\int _{-\pi}^{\pi}
\frac{|1-e^{i\lambda\tau}|^{2d}}{\lambda^{2d}}
\left\langle B_{1} ,f(\lambda )\right\rangle d\lambda =p\bigg\};
\end{multline*}
and
\begin{align*}
 {\md D_{V}^{U1}}&=\bigg\{g(\lambda )\bigg|V(\lambda )\leq g(\lambda
)\leq U(\lambda ), \frac{1}{2\pi } \int _{-\pi}^{\pi}
g(\lambda )d\lambda=Q\bigg\},
\\
  {\md D_{V}^{U2}} &=\bigg\{g(\lambda )\bigg|{\mathrm{Tr}}\, [V(\lambda
)]\leq {\mathrm{Tr}}\,[ g(\lambda )]\leq {\mathrm{Tr}}\, [U(\lambda )],
\frac{1}{2\pi } \int _{-\pi}^{\pi}
{\mathrm{Tr}}\,  [g(\lambda)]d\lambda  =q \bigg\},
\\
{\md D_{V}^{U3}} &=\bigg\{g(\lambda )\bigg|v_{kk} (\lambda )  \leq
g_{kk} (\lambda )\leq u_{kk} (\lambda ),
\frac{1}{2\pi} \int _{-\pi}^{\pi}
g_{kk} (\lambda
)d\lambda  =q_{k} , k=\overline{1,\infty}\bigg\},
\\
{\md D_{V}^{U4}} &=\bigg\{g(\lambda )\bigg|\left\langle B_{2}
,V(\lambda )\right\rangle \leq \left\langle B_{2},g(\lambda
)\right\rangle \leq \left\langle B_{2} ,U(\lambda)\right\rangle,
\frac{1}{2\pi }
\int _{-\pi}^{\pi}
\left\langle B_{2},g(\lambda)\right\rangle d\lambda  =q\bigg\}.
\end{align*}

\noindent
Here  $f_{1} ( \lambda )$, $V( \lambda )$, $U( \lambda )$ are known and fixed spectral densities, $W(\lambda)$ is an unknown spectral density, $p, p_k, q, q_k, k=\overline{1,\infty}$, are given numbers, $P, B_1, Q, B_2$ are given positive-definite Hermitian matrices.

The condition $0\in\partial\Delta_{\mathcal{D}}(f^0,g^0)$
implies the following equations  determining the least favourable spectral densities for these given sets of admissible spectral densities.

For the first set of admissible spectral densities $\md D_{\varepsilon}^{1}\times \md D_{V}^{U1}$, we have equations
\begin{multline}  \label{eq_5_1f_e_n_c}
\left(
{\me C}^{f0}_{\tau}(e^{i \lambda})
\right)
\left(
{\me C}^{f0}_{\tau}(e^{i \lambda})
\right)^{*}=
\\
=
\left(\frac{|1-e^{i\lambda\tau}|^{2d}}{\lambda^{2d}} (f^0(\lambda)+{\lambda}^{2d}g^0(\lambda))\right)
(\vec{\alpha}_f\cdot \vec{\alpha}^{*}+\Gamma(\lambda))
\times
\\
\times
\left(\frac{|1-e^{i\lambda\tau}|^{2d}}{\lambda^{2d}} (f^0(\lambda)+{\lambda}^{2d}g^0(\lambda))\right),
\end{multline}
\begin{multline} \label{eq_4_1g_e_n_c}
\left(
{\me C}^{g0}_{\tau}(e^{i \lambda})
\right)
\left(
{\me C}^{g0}_{\tau}(e^{i \lambda})
\right)^{*}=
\\
=
\left(\frac{|1-e^{i\lambda\tau}|^{2d}}{\lambda^{2d}}(f^0(\lambda)+{\lambda}^{2d}g^0(\lambda))\right)
(\vec{\beta}\cdot \vec{\beta}^{*}+\Gamma _{1} (\lambda )+\Gamma _{2} (\lambda ))
\times
\\
\times
\left(\frac{|1-e^{i\lambda\tau}|^{2d}}{\lambda^{2d}}(f^0(\lambda)+{\lambda}^{2d}g^0(\lambda))\right),
\end{multline}

\noindent where $\vec{\alpha}_f$ and $ \vec{\beta}$ are vectors of Lagrange multipliers, function $\Gamma(\lambda )\le 0$ and $\Gamma(\lambda )=0$ if $f_{0}(\lambda )>(1-\varepsilon )f_{1} (\lambda )$,
the matrix $\Gamma _{1} (\lambda )\le 0$ and $\Gamma _{1} (\lambda )=0$ if $g_{0}(\lambda )>V(\lambda ),$ the matrix  $
\Gamma _{2} (\lambda )\ge 0$ and $\Gamma _{2} (\lambda )=0$ if $g_{0}(\lambda )<U(\lambda ).$

For the second set of admissible spectral densities  $\md D_{\varepsilon}^{2}\times \md D_{V}^{U2}$, we have equations
\begin{equation} \label{eq_5_2f_e_n_c}
\left(
{\me C}^{f0}_{\tau}(e^{i \lambda})
\right)
\left(
{\me C}^{f0}_{\tau}(e^{i \lambda})
\right)^{*}
=
(\alpha_{f}^{2} +\gamma(\lambda ))
\left(\frac{|1-e^{i\lambda\tau}|^{2d}}{\lambda^{2d}} (f^0(\lambda)+{\lambda}^{2d}g^0(\lambda))\right)
^2,
\end{equation}
\begin{multline} \label{eq_4_2g_e_n_c}
\left(
{\me C}^{g0}_{\tau}(e^{i \lambda})
\right)
\left(
{\me C}^{g0}_{\tau}(e^{i \lambda})
\right)^{*}=
\\
=(\beta^{2} +\gamma _{1} (\lambda )+\gamma _{2} (\lambda )) \left(\frac{|1-e^{i\lambda\tau}|^{2d}}{\lambda^{2d}}(f^0(\lambda)+{\lambda}^{2d}g^0(\lambda))\right)^2,
\end{multline}

\noindent where  $\alpha_{f}^{2}$, $ \beta^{2}$ are Lagrange multipliers,  function $\gamma(\lambda )\le 0$ and $\gamma(\lambda )=0$ if ${\mathrm{Tr}}\,[f^{0} (\lambda )]>(1-\varepsilon ) {\mathrm{Tr}}\, [f_{1} (\lambda )]$,
 function $\gamma _{1} (\lambda )\le 0$ and $\gamma _{1} (\lambda )=0$ if ${\mathrm{Tr}}\,
[g^{0} (\lambda )]> {\mathrm{Tr}}\,  [V(\lambda )],$ the function $\gamma _{2} (\lambda )\ge 0$ and $\gamma _{2} (\lambda )=0$ if $ {\mathrm{Tr}}\,[g^{0}(\lambda )]< {\mathrm{Tr}}\, [ U(\lambda)].$

For the third set of admissible spectral densities $\md D_{\varepsilon}^{3}\times \md D_{V}^{U3}$, we have equation
\begin{multline}   \label{eq_5_3f_e_n_c}
\left(
{\me C}^{f0}_{\tau}(e^{i \lambda})
\right)
\left(
{\me C}^{f0}_{\tau}(e^{i \lambda})
\right)^{*}=
\\
=
\left(\frac{|1-e^{i\lambda\tau}|^{2d}}{\lambda^{2d}} (f^0(\lambda)+{\lambda}^{2d}g^0(\lambda))\right)
\left\{(\alpha_{fk}^{2} +\gamma_{k} (\lambda ))\delta _{kl} \right\}_{k,l=1}^{\infty}
\times
\\
\times
\left(\frac{|1-e^{i\lambda\tau}|^{2d}}{\lambda^{2d}} (f^0(\lambda)+{\lambda}^{2d}g^0(\lambda))\right),
\end{multline}
\begin{multline} \label{eq_4_3g_e_n_c}
\left(
{\me C}^{g0}_{\tau}(e^{i \lambda})
\right)
\left(
{\me C}^{g0}_{\tau}(e^{i \lambda})
\right)^{*}=
\\
=\left(\frac{|1-e^{i\lambda\tau}|^{2d}}{\lambda^{2d}}(f^0(\lambda)+{\lambda}^{2d}g^0(\lambda))\right)
\left\{(\beta_{k}^{2} +\gamma _{1k} (\lambda )+\gamma _{2k} (\lambda ))\delta _{kl}\right\}_{k,l=1}^{\infty}
\times
\\
\times
\left(\frac{|1-e^{i\lambda\tau}|^{2d}}{\lambda^{2d}}(f^0(\lambda)+{\lambda}^{2d}g^0(\lambda))\right),
\end{multline}

\noindent where  $\alpha_{fk}^{2}$,   $\beta_{k}^{2}$ are Lagrange multipliers,
 $\delta _{kl}$ are Kronecker symbols,
 functions $\gamma_{k}(\lambda )\le 0$ and $\gamma_{k} (\lambda )=0$ if $f_{kk}^{0}(\lambda )>(1-\varepsilon )f_{kk}^{1} (\lambda )$,
functions $\gamma _{1k} (\lambda )\le 0$ and $\gamma _{1k} (\lambda )=0$ if $g_{kk}^{0} (\lambda )>v_{kk} (\lambda ),$
functions $\gamma _{2k} (\lambda )\ge 0$ and $\gamma _{2k} (\lambda )=0$ if $g_{kk}^{0} (\lambda )<u_{kk} (\lambda).$

For the forth set of admissible spectral densities $\md D_{\varepsilon}^{4}\times \md D_{V}^{U4}$, we have equation
\begin{multline}   \label{eq_5_4f_e_n_c}
\left(
{\me C}^{f0}_{\tau}(e^{i \lambda})
\right)
\left(
{\me C}^{f0}_{\tau}(e^{i \lambda})
\right)^{*}=
\\
=
(\alpha_f^{2} +\gamma'(\lambda ))
\left(\frac{|1-e^{i\lambda\tau}|^{2d}}{\lambda^{2d}} (f^0(\lambda)+{\lambda}^{2d}g^0(\lambda))
\right)
 B_{1}^{\top}
 \times
\\
\times
\left(\frac{|1-e^{i\lambda\tau}|^{2d}}{\lambda^{2d}} (f^0(\lambda)+{\lambda}^{2d}g^0(\lambda))\right),
\end{multline}
\begin{multline} \label{eq_4_4g_e_n_c}
\left(
{\me C}^{g0}_{\tau}(e^{i \lambda})
\right)
\left(
{\me C}^{g0}_{\tau}(e^{i \lambda})
\right)^{*}=
\\
=
(\beta^{2} +\gamma'_{1}(\lambda )+\gamma'_{2}(\lambda ))
\left(\frac{|1-e^{i\lambda\tau}|^{2d}}{\lambda^{2d}}(f^0(\lambda)+{\lambda}^{2d}g^0(\lambda))\right)
B_{2}^{\top}
\times
\\
\times
\left(\frac{|1-e^{i\lambda\tau}|^{2d}}{\lambda^{2d}}(f^0(\lambda)+{\lambda}^{2d}g^0(\lambda))\right),
\end{multline}

\noindent where  $\alpha _{f}^{2}$, $ \beta^{2}$, are Lagrange multipliers, function $\gamma' ( \lambda )\le 0$ and $\gamma' ( \lambda )=0$ if $\langle B_{1} ,f^{0} ( \lambda ) \rangle>(1- \varepsilon ) \langle B_{1} ,f_{1} ( \lambda ) \rangle$,
functions $\gamma'_{1}( \lambda )\le 0$ and $\gamma'_{1} ( \lambda )=0$ if $\langle B_{2},g^{0} ( \lambda) \rangle > \langle B_{2},V( \lambda ) \rangle,$ functions $\gamma'_{2}( \lambda )\ge 0$ and $\gamma'_{2} ( \lambda )=0$ if $\langle
B_{2} ,g^{0} ( \lambda) \rangle < \langle B_{2} ,U( \lambda ) \rangle.$

The following theorem  holds true.

\begin{thm}
The least favorable spectral densities $f^{0}(\lambda)$, $g^{0}(\lambda)$ in classes
$\md D_{\varepsilon}^{k}\times \md D_{V}^{Uk},k=1,2,3,4$ for the optimal linear extrapolation of the functional  $A{\xi}$ from observations of the process $\xi(t)+\eta(t)$ at points  $t<0$ are determined by equations
\eqref{eq_5_1f_e_n_c}--\eqref{eq_4_1g_e_n_c},  \eqref{eq_5_2f_e_n_c}--\eqref{eq_4_2g_e_n_c}, \eqref{eq_5_3f_e_n_c}--\eqref{eq_4_3g_e_n_c}, \eqref{eq_5_4f_e_n_c}--\eqref{eq_4_4g_e_n_c},
respectively,  the minimality condition (\ref{umova11_e_n_c}), the
  constrained optimization problem (\ref{minimax1_e_n_c}) and restrictions  on densities from the corresponding classes
$ \md D_{\varepsilon}^{k}, \md D_{V}^{Uk},k=1,2,3,4$.  The minimax-robust spectral characteristic of the optimal estimate of the functional $A{\xi}$ is determined by  formula (\ref{spectr_A_e_n_c}).
\end{thm}

Let the spectral densities $f(\lambda)$ and $g(\lambda)$ admit the canonical factorizations   (\ref{fakt1_e_n_c}) and (\ref{fakt3_e_n_c}).
The   equations for the least favourable spectral densities are described below.

For the first set of admissible spectral densities $\md D_{\varepsilon}^{1}\times \md D_{V}^{U1}$, we have equation
\begin{align}  \label{eq_5_1f_fact_e_n_c}
\left(
{\me r^{0}_{\tau,f}(e^{i\lambda})}
\right)
\left(
{\me r^{0}_{\tau,f}(e^{i\lambda})}
\right)^{*}
&=
(\Theta_{\tau}(e^{-i\lambda}))^{\top}
(\vec{\alpha}_f\cdot \vec{\alpha}_f^{*}+\Gamma(\lambda))
\overline{\Theta_{\tau}(e^{-i\lambda})},
\\ \label{eq_4_1g_fact_e_n_c}
\left(
{\me r^{0}_{\tau,g}(e^{i\lambda})}
\right)
\left(
{\me r^{0}_{\tau,g}(e^{i\lambda})}
\right)^{*}
&=
(\Theta_{\tau}(e^{-i\lambda}))^{\top}
(\vec{\beta}\cdot \vec{\beta}^{*}+\Gamma _{1} (\lambda )+\Gamma _{2} (\lambda ))
\overline{\Theta_{\tau}(e^{-i\lambda})}
\end{align}
\noindent where $\vec{\alpha}_f$ and $ \vec{\beta}$ are vectors of Lagrange multipliers, matrix $\Gamma(\lambda )\le 0$ and $\Gamma(\lambda )=0$ if $f^{0}(\lambda )>(1-\varepsilon )f_{1} (\lambda )$,
 matrix $\Gamma _{1} (\lambda )\le 0$ and $\Gamma _{1} (\lambda )=0$ if $g^{0}(\lambda )>V(\lambda ),$  matrix  $
\Gamma _{2} (\lambda )\ge 0$ and $\Gamma _{2} (\lambda )=0$ if $g^{0}(\lambda )<U(\lambda ).$

For the second set of admissible spectral densities  $\md D_{\varepsilon}^{2}\times \md D_{V}^{U2}$, we have equations
\begin{align} \label{eq_5_2f_fact_e_n_c}
\left(
{\me r^{0}_{\tau,f}(e^{i\lambda})}
\right)
\left(
{\me r^{0}_{\tau,f}(e^{i\lambda})}
\right)^{*}
&=
(\alpha_f^{2} +\gamma(\lambda ))
(\Theta_{\tau}(e^{-i\lambda}))^{\top}\overline{\Theta_{\tau}(e^{-i\lambda})},
\\ \label{eq_4_2g_fact_e_n_c}
\left(
{\me r^{0}_{\tau,g}(e^{i\lambda})}
\right)
\left(
{\me r^{0}_{\tau,g}(e^{i\lambda})}
\right)^{*}
&=(\beta^{2} +\gamma _{1} (\lambda )+\gamma _{2} (\lambda )) (\Theta_{\tau}(e^{-i\lambda}))^{\top}\overline{\Theta_{\tau}(e^{-i\lambda})},
\end{align}

\noindent where  $\alpha _{f}^{2}$, $ \beta^{2}$ are Lagrange multipliers, function $\gamma(\lambda )\le 0$ and $\gamma(\lambda )=0$ if ${\mathrm{Tr}}\,[f^{0} (\lambda )]>(1-\varepsilon ) {\mathrm{Tr}}\, [f_{1} (\lambda )]$,
 function $\gamma _{1} (\lambda )\le 0$ and $\gamma _{1} (\lambda )=0$ if ${\mathrm{Tr}}\,
[g^{0} (\lambda )]> {\mathrm{Tr}}\,  [V(\lambda )],$ the function $\gamma _{2} (\lambda )\ge 0$ and $\gamma _{2} (\lambda )=0$ if $ {\mathrm{Tr}}\,[g^{0}(\lambda )]< {\mathrm{Tr}}\, [ U(\lambda)].$

For the third set of admissible spectral densities $\md D_{\varepsilon}^{3}\times \md D_{V}^{U3}$, we have equation
\begin{equation}   \label{eq_5_3f_fact_e_n_c}
\left(
{\me r^{0}_{\tau,f}(e^{i\lambda})}
\right)
\left(
{\me r^{0}_{\tau,f}(e^{i\lambda})}
\right)^{*}
=
(\Theta_{\tau}(e^{-i\lambda}))^{\top}
\left\{(\alpha_{fk}^{2} +\gamma_{k} (\lambda ))\delta _{kl} \right\}_{k,l=1}^{\infty}
\overline{\Theta_{\tau}(e^{-i\lambda})},
\end{equation}
\begin{multline} \label{eq_4_3g_fact_e_n_c}
\left(
{\me r^{0}_{\tau,g}(e^{i\lambda})}
\right)
\left(
{\me r^{0}_{\tau,g}(e^{i\lambda})}
\right)^{*}
=(\Theta_{\tau}(e^{-i\lambda}))^{\top}
\times\\
\times\left\{(\beta_{k}^{2} +\gamma _{1k} (\lambda )+\gamma _{2k} (\lambda ))\delta _{kl}\right\}_{k,l=1}^{\infty}
\overline{\Theta_{\tau}(e^{-i\lambda})},
\end{multline}

\noindent where  $\alpha _{fk}^{2}$,   $\beta_{k}^{2}$ are Lagrange multipliers,
 $\delta _{kl}$ are Kronecker symbols, functions $\gamma_{k}(\lambda )\le 0$ and $\gamma_{k}(\lambda )=0$ if $f_{kk}^{0}(\lambda )>(1-\varepsilon )f_{kk}^{1} (\lambda )$,
functions $\gamma _{1k} (\lambda )\le 0$ and $\gamma _{1k} (\lambda )=0$ if $g_{kk}^{0} (\lambda )>v_{kk} (\lambda ),$ functions $\gamma _{2k} (\lambda )\ge 0$ and $\gamma _{2k} (\lambda )=0$ if $g_{kk}^{0} (\lambda )<u_{kk} (\lambda).$

For the fourth set of admissible spectral densities $\md D_{\varepsilon}^{4}\times \md D_{V}^{U4}$, we have equation
\begin{align}   \label{eq_5_4f_fact_e_n_c}
\left(
{\me r^{0}_{\tau,f}(e^{i\lambda})}
\right)
\left(
{\me r^{0}_{\tau,f}(e^{i\lambda})}
\right)^{*}
&=
(\alpha_f^{2} +\gamma'(\lambda ))
(\Theta_{\tau}(e^{-i\lambda}))^{\top}
 B_{1}
 \overline{\Theta_{\tau}(e^{-i\lambda})},
\\ \label{eq_4_4g_fact_e_n_c}
\left(
{\me r^{0}_{\tau,g}(e^{i\lambda})}
\right)
\left(
{\me r^{0}_{\tau,g}(e^{i\lambda})}
\right)^{*}
&=
(\beta^{2} +\gamma'_{1}(\lambda )+\gamma'_{2}(\lambda ))
(\Theta_{\tau}(e^{-i\lambda}))^{\top}
B_{2}
\overline{\Theta_{\tau}(e^{-i\lambda})},
\end{align}

\noindent where $\alpha _{f}^{2}$, $ \beta^{2}$,   are Lagrange multipliers, function $\gamma' ( \lambda )\le 0$ and $\gamma' ( \lambda )=0$ if $\langle B_{1} ,f^{0} ( \lambda ) \rangle>(1- \varepsilon ) \langle B_{1} ,f_{1} ( \lambda ) \rangle$,
 functions $\gamma'_{1}( \lambda )\le 0$ and $\gamma'_{1} ( \lambda )=0$ if $\langle B_{2},g^{0} ( \lambda) \rangle > \langle B_{2},V( \lambda ) \rangle,$ functions $\gamma'_{2}( \lambda )\ge 0$ and $\gamma'_{2} ( \lambda )=0$ if $\langle
B_{2} ,g^{0} ( \lambda) \rangle < \langle B_{2} ,U( \lambda ) \rangle.$

The following theorems  hold   true.

\begin{thm}
The least favorable spectral densities $f^{0}(\lambda)$,  $g^{0}(\lambda)$  in the classes
$\md D_{\varepsilon}^{k}\times \md D_{V}^{Uk},k=1,2,3,4$ for the optimal linear prediction of the functional  $A{\xi}$ from observations of the process $\xi(t)+\eta(t)$ at points  $t<0$ by  canonical factorizations   (\ref{fakt1_e_n_c}) and (\ref{fakt3_e_n_c}),
 equations
\eqref{eq_5_1f_fact_e_n_c}--\eqref{eq_4_1g_fact_e_n_c},  \eqref{eq_5_2f_fact_e_n_c}--\eqref{eq_4_2g_fact_e_n_c}, \eqref{eq_5_3f_fact_e_n_c}--\eqref{eq_4_3g_fact_e_n_c}, \eqref{eq_5_4f_fact_e_n_c}--\eqref{eq_4_4g_fact_e_n_c},
respectively,
  constrained optimization problem (\ref{simple_minimax1_e_n_c}) and restrictions  on densities from the corresponding classes
$ \md D_{\varepsilon}^{k}, \md D_{V}^{Uk},k=1,2,3,4$.  The minimax-robust spectral characteristic of the optimal estimate of the functional $A{\xi}$ is determined by   formula (\ref{simple_spectr A_e_n_c}).
\end{thm}

\begin{thm}
 If the spectral density $g(\lambda)$ is known, the least favorable spectral density $f^{0}(\lambda)$ in the classes $\md D_{\varepsilon}^{k}$, $k=1,2,3,4$
 for the optimal linear prediction of the functional  $A{\xi}$ from observations of the process $\xi(t)+\eta(t)$ at points  $t<0$  is determined
by  canonical factorizations (\ref{fakt3_e_n_c}) and (\ref{fakt1_e_n_c}),
 equations
\eqref{eq_5_1f_fact_e_n_c},  \eqref{eq_5_2f_fact_e_n_c}, \eqref{eq_5_3f_fact_e_n_c}, \eqref{eq_5_4f_fact_e_n_c},
respectively,
 constrained optimization problem (\ref{simple_minimax3_e_n_c}) and restrictions  on density from the corresponding classes $\md D_{f\varepsilon}^{k}$, $k=1,2,3,4$.  The minimax-robust spectral characteristic of the optimal estimate of the functional $A{\xi}$ is determined by   formula (\ref{simple_spectr A_e_n_c}).
\end{thm}

\section*{Conclusions}

In this article, we dealt with  continuous time stochastic processes with periodically correlated $d$th increments.
These stochastic processes form a class of non-stationary stochastic processes
that combine  periodic structure of covariation functions of processes as well as integrating one.

We derived solutions of the problem of estimation of the linear functionals constructed from the  unobserved values of a continuous time stochastic process with periodically correlated $d$th increments.
Estimates are based on  observations of this process with periodically stationary noise at points $t<0$.
We obtained the estimates by representing the process under investigation as a vector-valued sequence with stationary increments.
Based on the  solutions for these type of sequences, we solved the corresponding problem for the considered class of continuous time stochastic processes.
The problem is investigated in the case of spectral certainty, where spectral densities of sequences are exactly known.
In this case  we propose an approach based on the Hilbert space projection method.
We derive formulas for calculating the spectral characteristics and the mean-square errors of the optimal estimates of the functionals.
In the case of spectral uncertainty where the spectral densities are not exactly known while, instead,
some sets of admissible spectral densities are specified, the minimax-robust method is applied.
We propose a representation of the mean square error in the form of a linear
functional in $L_1$ with respect to spectral densities, which allows
us to solve the corresponding constrained optimization problem and
describe the minimax-robust estimates of the functionals. Formulas
that determine the least favorable spectral densities and minimax-robust spectral characteristics of the optimal linear estimates of
the functionals are derived  for a wide list of specific classes
of admissible spectral densities.

\section*{Appendix}
\emph{\textbf{Proof of Lemma \ref{predst H_cont}.}}

 Define
\[
b^{\tau}_j(u)=b^{\tau}(u+jT),\quad \zeta^{(d)}_j(u)=\zeta^{(d)}_j(u+jT,\tau T),\,\, u\in [0,T),
\]
and
\[
a_j(u)=a(u+jT),\quad \eta_j(u)=\eta_j(u+jT),\,\, u\in [0,T).
\]

Making use of  decomposition (\ref{zeta}) for the increment sequence  $\{\zeta^{(d)}_j,j\in\mathbb Z\}$  and the solution of equation
\begin{equation} \label{peretv}
(-1)^k\left[\frac{k}{2}\right]+(-1)^m\left[\frac{m}{2}\right]=0
\end{equation}
of two variables  $(k,m)$, which is given by pairs  $(1,1)$, $(2l+1,2l)$ and
$(2l,2l+1)$ for $l=2,3,\dots$, rewrite the functional $B\zeta$ as
\cite{MoklyachukGolichenko2016}
\begin{align*}
    B\zeta&=\int_0^{\infty} b^{\tau}(t)\zeta^{(d)}(t,\tau T)dt
= \sum_{j=0}^{\infty}\int_{0}^{T}
b^{\tau}_j(u)\zeta^{(d)}_j(u)
du
\\
&=\sum_{j=0}^{\infty}\frac{1}{T}\int_{0}^{T}\left(\sum_{k=1}^\infty b^{\tau}_{kj}
e^{2\pi i\{(-1)^k\left[\frac{k}{2}\right]\}u/T} \right)
\left(\sum_{m=1}^\infty \zeta^{(d)}_{mj} e^{2\pi
i\{(-1)^m\left[\frac{m}{2}\right]\}u/T} \right)du
\\
&=\sum_{j=0}^{\infty}\sum_{k=1}^\infty \sum_{m=1}^\infty b^{\tau}_{kj} \zeta^{(d)}_{mj}
\frac{1}{T}\int_{0}^{T}e^{2\pi
i\left\{(-1)^k\left[\frac{k}{2}\right]+(-1)^m\left[\frac{m}{2}\right]\right\}u/T}
du
\\
&=\sum_{j=0}^{\infty}\sum_{k=1}^{\infty} b^{\tau}_{kj}\zeta^{(d)}_{kj}=
 \sum_{j=0}^{\infty}
{(\vec{b}^{\tau}_j)}^{\top}\vec{\zeta}^{(d)}_j
\\
&=B\vec \zeta.
\end{align*}

The representation of the  functional $V\eta$ is obtained in the same way:
\begin{align*}
    V\eta&=\int_0^{\infty} a(t)\eta(t)dt
= \sum_{j=0}^{\infty}\int_{0}^{T}
a_j(u)\eta_j(u)
du
\\
&=\sum_{j=0}^{\infty}\sum_{k=1}^{\infty} a_{kj}\eta_{kj}=
 \sum_{j=0}^{\infty}
{(\vec{a}_j)}^{\top}\vec{\eta}_j
\\
&=V\vec \eta.
\end{align*}
From Lemma \ref{predst A_cont} we obtain
\[
b^{\tau}_j(u)=\sum_{l=0}^{\infty}a( u+jT +\tau T l)d(l)= D^{\tau T} a(u),\,
u\in[0;T), \, j=0,1,\ldots,
\]
and
\[
b^{\tau}_{kj}=\sum_{l=0}^{\infty}a_{kj+\tau l}d(l),\,
 j=0,1,\ldots,
\]
which finalizes the proof of Lemma \ref{predst H_cont}.

\

\emph{\textbf{Proof of Theorem \ref{thm1_e_n_c}.}}

A projection $\widehat{H}\vec\xi$ of the element $H\vec\xi$ on the
subspace $H^{0-}(\xi^{(d)}_{\tau}+\eta^{(d)}_{\tau})$  is
characterized by two conditions:

1) $ \widehat{H}\vec\xi\in H^{0-}(\xi^{(d)}_{\tau}+\eta^{(d)}_{\tau})$;

2) $(H\vec\xi-\widehat{H}\vec\xi)
\perp
H^{0-}(\xi^{(d)}_{\tau}+\eta^{(d)}_{\tau})$.

 The second condition implies the following relation which holds true for all $j\leq-1$ and for all $k\geq1$
 \begin{multline*}
\int_{-\pi}^{\pi}
\bigg((\vec{B}_{\tau}(e^{i\lambda}))^{\top}\frac{(1-e^{-i\lambda\tau})^d}{(i\lambda)^d}(f(\lambda)+\lambda^{2d}g(\lambda))-
(\vec{h}_{\tau}(\lambda))^{\top}
(f(\lambda)+\lambda^{2d}g(\lambda))
-\bigg.
\\
\bigg.-
(\vec{A}(e^{i\lambda}))^{\top}g(\lambda)(-i\lambda)^d\bigg)\vec{\delta}_{k}
\frac{(1-e^{i\lambda\tau})^d}{(-i\lambda)^{d}}e^{-i\lambda j}d\lambda=0.
 \end{multline*}
From these relations, we conclude that the spectral characteristic $\vec{h}_{\tau}(\lambda)$ of the estimate
$\widehat{H}\vec\xi$ allow a representation in the form  (\ref{spectr_A_e_n_c})
where
\[
\vec{C}_{\tau}(e^{i \lambda})=\sum_{j=0}^{\infty}\vec{c}_j^{\tau}e^{i\lambda j},
\]
and $\vec{c}_j^{\tau}=\{c_{kj}^{\tau}\}_{k=1}^{\infty}$, $j=0,1,2,\dots$,  are unknown coefficients to be found.

Condition 1)  implies
${(i\lambda)^d }{(1-e^{-i\lambda \tau})^{-d}}\vec{h}_{\tau}(\lambda)\in L_2^{0-},
$
and thus,
\begin{multline*}
\int_{-\pi}^{\pi} \biggl[(\vec{B}_{\tau}(e^{i\lambda}))^{\top}-
(\vec{A}_{\tau}(e^{i\lambda}))^{\top}g(\lambda)
\frac{\lambda^{2d}}{|1-e^{i\lambda \tau}|^{2d}}
(f(\lambda)+{\lambda}^{2d}g(\lambda))^{-1}
\\
-
(\vec{C}_{\tau}(e^{i\lambda}))^{\top}\frac{\lambda^{2d}}{|1-e^{i\lambda \tau}|^{2d}}(f(\lambda)+{\lambda}^{2d}g(\lambda))^{-1}\biggr]e^{-i\lambda j}d\lambda=\vec 0,\quad  j\geq 0,
\label{eq_C_f_n_c}
\end{multline*}
which can be presented as a system of linear equations
 \be \label{linear equations1_e_n_c}
  \vec{b}_l^{\tau}-\sum_{ j=0}^{\infty}T^{\tau}_{l,j}\vec{a}_j^{\tau}
    =\sum_{j=0}^{\infty}P_{l,j}^{\tau}\vec{c}_j^{\tau},\quad l\geq0,
    \ee
  determining the unknown coefficients ${\vec c}^{\tau}_j$, $j\geq0$.

Rewrite the  system of equations (\ref{linear equations1_e_n_c}) in the matrix form
\[D^{\tau}\me a-\me T_{\tau}\me a^{\tau}=\me P_{\tau}\me c^{\tau},\]
 where
 \[
\me c^{\tau}=((\vec{c}^{\tau}_0)^{\top},(\vec{c}^{\tau}_1)^{\top},(\vec{c}^{\tau}_2)^{\top}, \ldots)^{\top}.
 \]
Consequently, the unknown coefficients $\vec{c}^{\tau}_j$, $j\geq0$,  determining
the spectral characteristic $\vec{h}_{\tau}(\lambda)$ are as follows
\begin{equation*}\label{meq-e_n_c}
\vec{c}^{\tau}_j=(\me P_{\tau}^{-1}D^{\tau}\me a-\me P_{\tau}^{-1}\me T_{\tau}\me a^{\tau})_j,\quad j\geq 0,
 \end{equation*}
where $(\me P_{\tau}^{-1}D^{\tau}\me a-\me
P_{\tau}^{-1}\me T_{\tau}\me a^{\tau})_j$, $j\geq 0$, is the
$j$th  infinite dimension vector entry of the vector  $\me P_{\tau}^{-1}D^{\tau}\me a-\me P_{\tau}^{-1}\me T_{\tau}\me a^{\tau}$.
The existence of the inverse matrix $(\me P^{\tau})^{-1}$ is justified  in \cite{Luz_Mokl_book} under condition (\ref{umova11_e_n_c}).
Thus, the function $\vec{C}_{\tau}(e^{i \lambda})$ is calculated as
\[
\vec{C}_{\tau}(e^{i \lambda})=\sum_{j=0}^{\infty}
(\me P_{\tau}^{-1}D^{\tau}\me a-\me P_{\tau}^{-1}\me T_{\tau}\me a^{\tau})_j
e^{i\lambda j}
\]
and the spectral characteristic $\vec{h}_{\tau}(\lambda)$ is calculated by the formula
\begin{align*}
\notag (\vec{h}_{\tau}(\lambda))^{\top}&=(\vec{B}_{\tau}(e^{i\lambda}))^{\top}
\frac{(1-e^{-i\lambda\tau})^d}{(i\lambda)^d}
\\
\notag&\quad-
(\vec{A}_{\tau}(e^{i\lambda}))^{\top}g(\lambda)\frac{(-i\lambda)^d}{(1-e^{i\lambda \tau})^d}
(f(\lambda)+{\lambda}^{2d}g(\lambda))^{-1}
\\
\notag&\quad-
\frac{(-i\lambda)^d}{(1-e^{i\lambda \tau})^d}
\left(
\sum_{k=0}^{\infty}(\me P_{\tau}^{-1}D^{\tau}\me a-\me P_{\tau}^{-1}\me T_{\tau}\me a^{\tau})_k e^{i\lambda j}
\right)^{\top}
\\&\quad\quad\times
(f(\lambda)+{\lambda}^{2d}g(\lambda))^{-1}.
\end{align*}
The value of the mean square error of the estimate $\widehat{A}\vec\xi$ is calculated by the formula
\begin{align*}
\notag& \Delta(f,g;\widehat{A}\vec\xi)=\Delta(f,g;\widehat{H}\vec\xi)= \mt E|H\vec\xi-\widehat{H}\vec\xi|^2=
\\
\notag&=
\frac{1}{2\pi}\int_{-\pi}^{\pi}
\frac{\lambda^{2d}}{|1-e^{i\lambda\tau}|^{2d}}\left[(\vec{A}_{\tau}(e^{i\lambda}))^{\top}g(\lambda) +
(\vec{C}_{\tau}(e^{i \lambda})
)^{\top}
\right]
\\\notag&\quad\quad
\times
(f(\lambda)+{\lambda}^{2d}g(\lambda))^{-1}\, f(\lambda)\, (f(\lambda)+{\lambda}^{2d}g(\lambda))^{-1}
\\\notag&\quad\quad
\times
\left[g(\lambda)\overline{\vec{A}_{\tau}(e^{i\lambda})} +
\overline{\vec{C}_{\tau}(e^{i \lambda})}
\right]
d\lambda
\\\notag&\quad
+\frac{1}{2\pi}\int_{-\pi}^{\pi}\frac{\lambda^{4d}}{|1-e^{i\lambda\tau}|^{4d}}
\left[\frac{|1-e^{i\lambda \tau}|^{2d}}{\lambda^{2d}}(\vec{A}(e^{i\lambda}))^{\top}f(\lambda)
-{(1-e^{-i\lambda \tau})^d}
(\vec{C}_{\tau}(e^{i \lambda}))^{\top}
\right]
\\
\notag&\quad\quad
\times
(f(\lambda)+{\lambda}^{2d}g(\lambda))^{-1}\,g(\lambda)\,(f(\lambda)+{\lambda}^{2d}g(\lambda))^{-1}
\\
\notag&\quad\quad
\times
\left[\frac{|1-e^{i\lambda \tau}|^{2d}}{\lambda^{2d}}f(\lambda)\overline{\vec{A}(e^{i\lambda})} -
{(1-e^{i\lambda \tau})^d}\overline{\vec{C}_{\tau}(e^{i \lambda}}
)
\right]
d\lambda
\\&=
\langle D^{\tau}\me a- \me T_{\tau}\me
 a^{\tau},\me P_{\tau}^{-1}D^{\tau}\me a-\me P_{\tau}^{-1}\me T_{\tau}\me a^{\tau}\rangle+\langle\me Q\me a,\me
 a\rangle,
\end{align*}
which finalizes the proof of Theorem \ref{thm1_e_n_c}.


\begin{thebibliography}{999}


\bibitem{Basawa}
\newblock I.V. Basawa, R. Lund, and Q. Shao,
\newblock \emph{First-order seasonal autoregressive processes with periodically varying parameters},
\newblock Statistics \& Probability Letters, vol. 67, no. 4, pp. 299--306, 2004.


\bibitem{DubovetskaMoklyachuk2013}
\newblock I. I. Dubovets'ka, and M. P. Moklyachuk,
\newblock \emph{Minimax estimation problem for periodically correlated stochastic processes},
\newblock Journal of Mathematics and System Science, vol. 3, no.1, pp. 26--30, 2013.

\bibitem{Dubovetska6}
\newblock I. I. Dubovets'ka, and M. P. Moklyachuk,
\newblock \emph{Extrapolation of periodically correlated processes from observations with noise},
\newblock Theory of Probability and Mathematical Statistics, vol. 88, pp. 67--83, 2014.


\bibitem{Dudek-Hurd}
\newblock A. Dudek, H. Hurd, and W. Wojtowicz,
\newblock \emph{PARMA methods based on Fourier representation of periodic coefficients},
\newblock Wiley Interdisciplinary Reviews: Computational Statistics, vol. 8, no. 3, pp. 130--149, 2016.

\bibitem{Franke1984}
\newblock J. Franke, and  H. V. Poor,
\newblock \emph{Minimax-robust filtering and finite-length robust predictors},
\newblock  In: Robust and Nonlinear Time
Series Analysis, Lecture Notes in Statistics, Springer-Verlag, No.26, pp. 87--126, 1984.


\bibitem{Gihman_Skorohod}
\newblock I. I. Gikhman, and A. V. Skorokhod,
\newblock Introduction to the theory of random processes.
\newblock Moskva: Gosudarstv. Izdat. Fiz.-Mat. Lit. 1965.


\bibitem{Glad1963}
 \newblock E. G. Gladyshev,
 \newblock \emph{Periodically and almost-periodically correlated random processes with continuous time parameter},
 \newblock Theory Probab. Appl. vol. 8, pp. 173--177, 1963.




\bibitem{Grenander}
\newblock U. Grenander,
\newblock \emph{A prediction problem in game theory},
\newblock Arkiv f\"or Matematik, vol. 3, pp. 371--379, 1957.


\bibitem{Hosoya}
 \newblock Y. Hosoya,
\newblock \emph{Robust linear extrapolations of second-order stationary processes},
 \newblock Annals of Probability, vol. 6, no. 4, pp. 574--584, 1978.


\bibitem{Kallianpur}
\newblock G. Kallianpur, and V. Mandrekar,:
\newblock  \emph{Spectral theory of stationary H-valued processes},
\newblock J. Multivariate Analysis, vol. 1, pp. 1--16,  1971.




\bibitem{Kassam_Poor1985}
 \newblock S. A. Kassam, and H. V. Poor,
\newblock \emph{Robust techniques for signal processing: A survey},
 \newblock Proceedings of the IEEE No.73, pp. 433-481, 1985.


\bibitem{Karhunen}
\newblock K. Karhunen,
\newblock \emph{Uber lineare Methoden in der Wahrscheinlichkeitsrechnung},
\newblock Annales Academiae Scientiarum Fennicae. Ser. A I, no. 37, 1947.


\bibitem{Kolmogorov}
\newblock A. N. Kolmogorov,
\newblock \emph{Selected works by A. N. Kolmogorov. Vol. II: Probability theory and mathematical statistics. Ed. by A.
N. Shiryayev. Mathematics and Its Applications. Soviet Series. 26. Dordrecht etc.}
\newblock Kluwer Academic Publishers, 1992.


\bibitem{Kozak_Mokl}
\newblock P. S. Kozak  and M. P. Moklyachuk,
\newblock \emph{Estimates of functionals constructed from random sequences with periodically stationary increments},
\newblock Theory Probability and Mathematical Statistics, vol. 97, pp. 85--98, 2018.


\bibitem{Luz_Mokl_extra}
\newblock M. Luz, and M. Moklyachuk,
\newblock \emph{Minimax-robust prediction problem for stochastic sequences with stationary increments and cointegrated sequences},
\newblock Statistics, Optimization and Information Compututing, vol. 3, no. 2, pp. 160--188, 2015.

\bibitem{Luz_Mokl_book}
\newblock M. Luz, and M. Moklyachuk,
\newblock \emph{Estimation of Stochastic Processes with Stationary Increments and Cointegrated Sequences},
\newblock London: ISTE; Hoboken, NJ: John Wiley \& Sons, 282 p., 2019.

\bibitem{Luz_Mokl_extra_GMI}
\newblock M. Luz, and M. Moklyachuk,
\newblock \emph{Minimax-robust forecasting of sequences with periodically stationary long memory multiple seasonal increments},
\newblock Statistics, Optimization and Information Computing, vol. 8, no. 3, pp. 684--721, 2020.

\bibitem{Luz_Mokl_extra_noise_PCI}
\newblock M. Luz, and M. Moklyachuk,
\newblock \emph{Minimax  prediction of sequences with periodically stationary increments observed with noise and cointegrated sequences},
\newblock In: M. Moklyachuk (ed.) Stochastic Processes: Fundamentals and Emerging Applications. Nova Science Publishers,  New York, pp. 189--247, 2023.

\bibitem{Luz_Mokl_extra_cont_PCI}
\newblock M. Luz, and M. Moklyachuk,
\newblock \emph{Estimation problem for continuous time stochastic processes with periodically correlated increments},
\newblock Statistics, Optimization and Information Computing, vol.11, no. 4, pp. 811 -- 828, 2023.


\bibitem {Moklyachuk:1981}
\newblock M. P. Moklyachuk,
\newblock \emph{Estimation of linear functionals of stationary stochastic processes and a two-person zero-sum game},
\newblock Stanford University Technical Report, no. 169, pp. 1--82, 1981.

\bibitem {Moklyachuk1981tv}
 \newblock M. P. Moklyachuk,
 \newblock \emph{On a problem of game theory and the extrapolation of stochastic processes with values in a Hilbert space},
\newblock Theory Probability and Mathematical Statistics, vol. 24, pp. 115--124, 1981.

\bibitem {Moklyachuk1982}
 \newblock M. P. Moklyachuk,
\newblock \emph{On an antagonistic game and prediction of stationary random sequences in a Hilbert space},
\newblock Theory Probability and Mathematical Statistics, vol. 25, pp. 107--113, 1982.



\bibitem{Moklyachuk:2008}
\newblock M. P. Moklyachuk,
\newblock \emph{Robust estimations of functionals of stochastic processes},
\newblock {Ky\"{\i}v: Vydavnychyj Tsentr ``Ky\"{\i}vs'ky\u{\i} Universytet''}, 320 p., 2008.  (in Ukrainian)


\bibitem{Moklyachuk2015}
\newblock M. P. Moklyachuk,
\newblock \emph{Minimax-robust estimation problems for stationary stochastic sequences},
\newblock Statistics, Optimization and Information Computing, vol. 3, no. 4, pp. 348--419, 2015.

\bibitem{MoklyachukGolichenko2016}
\newblock M. P. Moklyachuk, and I. I. Golichenko,
\newblock \emph{Periodically correlated processes estimates},
\newblock Saarbr\"ucken: LAP Lambert Academic Publishing, 308 p., 2016.



\bibitem{Mokl_Mas_book}
\newblock M.P. Moklyachuk, and A.Yu. Masyutka,
\newblock \emph{Minimax-robust estimation technique: For  stationary stochastic processes},
\newblock LAP Lambert Academic Publishing,  296 p. 2012.


\bibitem{Sidei_book}
\newblock M. Moklyachuk, M. Sidei, and O. Masyutka,
\newblock \emph{Estimation of stochastic processes with missing observations},
\newblock Mathematics Research Developments. New York, NY: Nova Science Publishers, 336 p., 2019.

\bibitem{Napolitano}
\newblock A. Napolitano,
\newblock \emph{Cyclostationarity: New trends and applications},
\newblock Signal Processing, vol. 120, pp. 385--408, 2016.


\bibitem{Pinsker}
\newblock M. S. Pinsker, and A. M. Yaglom,
\newblock \emph{On linear extrapolaion of random processes with $n$th stationary incremens},
\newblock Doklady Akademii Nauk SSSR, n. Ser. vol. 94, pp. 385--388, 1954.



\bibitem{Reisen2018}
\newblock V. A. Reisen,  E. Z. Monte,  G. C. Franco,  A. M. Sgrancio,  F. A. F. Molinares,  P. Bondond, F. A. Ziegelmann, and B. Abraham,
\newblock \emph{Robust estimation of fractional seasonal processes: Modeling and forecasting daily average SO2 concentrations},
\newblock Mathematics and Computers in Simulation, vol. 146, pp. 27--43, 2018.



\bibitem{Rockafellar}
\newblock R. T. Rockafellar,
\newblock \emph{Convex Analysis},
\newblock Princeton Landmarks in Mathematics. Princeton, NJ: Princeton University Press, 451 p., 1997.

\bibitem{Solci}
\newblock C. C. Solci, V. A. Reisen, A. J. Q. Sarnaglia, and  P. Bondon,
\newblock \emph{Empirical study of robust estimation methods for PAR models with application to the air quality area},
\newblock Communication in Statistics - Theory and Methods,  vol. 48, no. 1, pp. 152--168, 2020.


\bibitem{VastPoor1983}
\newblock S. K. Vastola, and H. V. Poor,
\newblock \emph{An analysis of the effects of spectral uncertainty on Wiener filtering},
\newblock Automatica vol.28, 289-293, 1983.


\bibitem{bookWien}
\newblock N. Wiener,
\newblock \emph{Extrapolation, interpolation, and smoothing of stationary time series. With engineering applications}.
\newblock  Cambridge, Mass.: The M. I. T. Press, Massachusetts Institute of Technology, 163 p. 1966.

\bibitem{Yaglom:1955}
\newblock  A. M. Yaglom,
 \newblock \emph{Correlation theory of stationary and related random processes with stationary $n$th increments},
\newblock Mat. Sbornik, vol. 37, no. 1, pp. 141--196, 1955.

\bibitem{Yaglom}
\newblock A. M.   Yaglom,
\newblock \emph{Correlation theory of stationary and related random functions. Vol. 1: Basic results;  Vol. 2: Suplementary notes and references},
\newblock Springer Series in Statistics, Springer-Verlag, New York etc., 1987.




\end{thebibliography}
\end{document}